	\newcommand{\bfig}{\begin{center}\vskip 0.5em}
	\newcommand{\efig}{\end{center}\vskip 0.5em}

 	\newcommand{\BOX}{\hbox {$\sqcap$ \kern -1em $\sqcup$}}

	\newcommand{\tensor}{\otimes}
	\newcommand\C{{\bf C}}
	\newcommand\R{{\bf R}}
	\newcommand{\be}{\begin{equation}}
     \newcommand{\ee}{\end{equation}}
     \newcommand{\ba}{\begin{eqnarray}}
     \newcommand{\ea}{\end{eqnarray}}
     \newcommand{\ban}{\begin{eqnarray*}}
     \newcommand{\ean}{\end{eqnarray*}}

	\newcommand{\Cob}{{\rm Cob}}
	\newcommand{\Vect}{{\rm Vect}}
	\newcommand{\Hilb}{{\rm Hilb}}
     \newcommand{\maps}{\colon}
	\newcommand{\tr}{{\rm tr}}

	\renewcommand{\hom}{{\rm hom}}
	\newcommand{\cat}{{\rm Cat}}
	\newcommand{\set}{{\rm Set}}
     \newcommand{\Cat}{{\rm Cat}}
	\newcommand{\Set}{{\rm Set}}
	
	\newcommand{\hot}{{\rm Hot}}
	\newcommand{\gpd}{{\rm Gpd}}

	\newcommand{\eps}{\epsilon}

	\newcommand{\doublearrow}{\Rightarrow}

	\documentclass{article}
        \usepackage{diagram}
        \usepackage[all, knot]{xy}
      \xyoption{arc}

\begin{document}

	\begin{center}
	{\bf  Higher-Dimensional Algebra and \\
         Topological Quantum Field Theory \\}
	\vspace{0.5cm}
	{\em John\ C.\ Baez and James Dolan\\}
	\vspace{0.3cm}
	{\small Department of Mathematics \\
	University of California\\
        Riverside CA 92521\\}
	\vspace{0.3cm}
	{\small February 5, 1995 \\}
	\vspace{0.3cm}
	\end{center}

\begin{abstract}

The study of topological quantum field theories  increasingly
relies upon concepts from higher-dimensional algebra such as
$n$-categories and $n$-vector spaces.   We review progress
towards a definition of $n$-category suited for this purpose, and
outline a program in which $n$-dimensional TQFTs are to be
described as $n$-category representations.   First we describe a
`suspension' operation on  $n$-categories, and hypothesize that
the $k$-fold suspension of a weak $n$-category stabilizes for $k
\ge n + 2$.   We give evidence for this hypothesis and describe
its relation to stable homotopy theory.  We then propose a
description of $n$-dimensional unitary extended TQFTs as weak
$n$-functors from the `free stable weak $n$-category with duals
on one object' to the $n$-category of `$n$-Hilbert spaces'.  We
conclude by  describing $n$-categorical generalizations of
deformation  quantization and the quantum double construction.

\end{abstract}

\section{Introduction}

One important lesson we have learned from topological quantum
field theory is that describing dynamics using {\it group}
representations is only a special case of describing it using
{\it category} representations.  In fact, examining the structure
of the known $n$-dimensional topological quantum field theories,
it appears that many are representations of some sort of {\it
$n$-category}.  While a category is a structure with objects and
morphisms between these objects, in an $n$-category there are
also `morphisms between morphisms' or 2-morphisms, `morphisms
between 2-morphisms' or 3-morphisms, and so on, up to
$n$-morphisms.  In the theory of manifolds, $k$-morphisms
correspond to manifolds (with boundary, corners, etc.) of
dimension $k$.

The theory of $n$-categories is one of several related approaches to
describing topology in purely algebraic terms.  Taken together, these
constitute a subject known as `higher-dimensional algebra'.  While
the basic insights of this subject are simple and beautiful, it is far
from reaching its final form.  The aim of this paper is to paint,
in rather broad strokes, a picture of some patterns that are
becoming apparent.  In many cases these are only well-understood in low
dimensions.  The topology of higher dimensions is known to be very
different, so uncautious extrapolation is risky.  Nonetheless, the
patterns we see so far should serve as a useful guide to research, if
only to goad us to a better understanding --- and in particular, an {\it
algebraic} understanding --- of the relation between topological
quantum field theory and the traditional techniques of algebraic
topology, which work in arbitrarily high dimensions.

In the rest of this section we review some of the physics issues
addressed by topological quantum field theories (TQFTs), and recall the
mathematical definition of a TQFT.  In Section 2, we give a sketchy
overview of why $n$-categories should serve as a natural framework for
an algebraic approach to TQFTs. In Section 3 we recall the definition of
a `strict $n$-category', and in Section 4 we begin to explain why a
weakening of this notion is crucial --- though only well-understood for
$n \le 3$.  In Section 5 we describe how to `suspend' an $n$-category,
and argue that the process of iterated suspension of an $n$-category
should stabilize after $n+2$ times.  In Section 6 we recall the roots of
higher-dimensional algebra in algebraic topology, and explain how our
notion of suspension relates to that in homotopy theory.  In Section 7
we argue that an algebraic framework for smooth manifold theory of
dimension $\le n$ is given by the `free stable $n$-category with duals
on one object'.  In Section 8 we then argue that a unitary extended
$n$-dimensional TQFT is a representation of this $n$-category in the
$n$-category of `$n$-Hilbert spaces'.  In Section 9 we conclude with
some examples illustrating the role of `quantization' in
higher-dimensional algebra.  Since we wish to assume a minimum of
familiarity with the subject, we take a rather expository approach, and
frequently refer the reader to review articles and books rather than
original papers.  For this reason, our bibliography is by no means
complete.

Why are TQFTs interesting as {\it physics}?  One reason is that
they possess certain features one expects of a quantum theory of
gravity.  This is not to say that quantum gravity is, or will be,
a TQFT.  Nonetheless, to understand the significance of TQFTs, it
is useful to recall some ideas from quantum gravity.

The `general covariance' or `diffeomorphism-invariance' of general
relativity is often regarded as a crucial feature to be preserved, if
possible, by any quantum theory accomodating gravity.  However, it has
traditionally been a bit unclear what we really are asking for when we
say we want a theory to be generally covariant.  In Einstein's original
work on general relativity \cite{Einstein}, he emphasized that the
equations should preserve their form under arbitrary coordinate
transformations: {\sl ``The general laws of nature are to be expressed
by equations which hold good for all systems of coordinates, that is,
are covariant with respect to any substitutions whatever (generally
covariant).''}

Later, it became clear that this definition of general covariance
is vacuous without some restriction on how the quantities
involved in the laws transform under coordinate transformations.
The importance of {\it tensors} and their transformation rules is
often stressed.  In fact, requiring all the quantities in a field
theory to be tensors misses the point in at least two ways.  On
the one hand, modern field theories typically make use of
nontensorial objects such as spinor fields, connections, and
other bundle sections. (Note that connections can be regarded as
bundle sections with the aid of jet bundles.)  On the other hand,
special-relativistic classical field theories on Minkowski space,
which we normally do not think of as `generally covariant,' can
nonetheless be reformulated as diffeomorphism-invariant, purely
tensorial equations by coupling them to Einstein gravity, setting
the gravitational constant equal to zero, and then adding an
extra equation saying that the metric is flat!

Reflection along these lines led to the recognition that the key
feature distinguishing general relativity from previous theories
is actually its lack of a fixed prior geometry.  In the words of
Misner, Thorne and Wheeler \cite{MTW}, {\sl ``By `prior geometry'
one means any aspect of the geometry of spacetime that is fixed
immutably, i.e., that cannot be changed by changing the
distribution of gravitating sources.''} The best-known example of
prior geometry is the Minkowski metric in special relativity.  In
special-relativistic quantum field theory, dynamics is described
in terms of representations of the symmetry group of this metric,
that is, the Poincar\'e group. Note that even in the
diffeomorphism-invariant reformulation of special-relativistic
field theory described above, there is, for each solution of the
equations of motion, a canonical choice of a subgroup of the
diffeomorphism group isomorphic to the Poincar\'e group: the
isometry group of the metric.

More generally, while {\it sufficient} conditions for a theory to
lack prior geometry are difficult to state, it certainly seems
{\it necessary} that there be no canonical way to choose, for
each state of the theory, a subgroup of the spacetime
diffeomorphism group isomorphic to some fixed nontrivial Lie
group.  For, by the Erlangen philosophy that a geometry is known
by its group of symmetries, such a subgroup would indicate the
presence of prior geometry independent of the state.  This, in
turn, suggests that in theories without prior geometry, dynamics
will not be described in terms of the representations of a Lie
group of spacetime symmetries.

The question is then, how {\it is} dynamics described in such
theories?  This is known as the `problem of time'.   In canonical
quantum gravity on a spacetime of the form $\R \times S$, it is
manifested by the fact that physical states should satisfy the
Wheeler-DeWitt equation
\[           H \psi = 0 ,\]
expressing their invariance under the diffeomorphism group.
This gives rise to the so-called `inner product problem', namely:
what is the correct inner product on the space of states?   In
quantum theories with prior geometry, the inner product is chosen
so that spacetime symmetry group acts unitarily on the space of
states, but in a formalism where states are
diffeomorphism-invariant the inner product must be determined in
some other way.

Many different approaches have been proposed to both these
problems \cite{Ashtekar,Isham}, but here we only consider one,
namely Atiyah's definition of a TQFT \cite{Atiyah}.  This is in
fact a very radical approach!  First, rather than attempting to
describe the dynamics of fields on a single spacetime manifold, a
TQFT describes the dynamics of fields in terms of a category in
which the objects are $(n-1)$-dimensional manifolds representing
possible choices of `space', and morphisms are $n$-dimensional
cobordisms representing choices of `spacetime'.  In fact, a TQFT
is a kind of a representation of this category, assigning a
vector space of states to each object and a linear operator to
each morphism.  Second, in any unitary TQFT satisfying
a certain nondegeneracy condition, the structure of
this category automatically determines the inner product in the
space of states for any $(n-1)$-manifold.

To fully appreciate these ideas one must recognize that a group is
just a very special sort of category.  Recall that a category consists
of a collection of objects, and for any objects $x$ and $y$, a set
$\hom(x,y)$ of morphisms from $x$ to $y$.  (Technically, we are
considering only locally small categories.)  If $f\in\hom(x,y)$, we
write $f\maps x\to y$.  Morphisms $g\maps x \to y$ and $f \maps y \to
z$ can be composed to obtain a morphism $fg
\maps x \to z$, and composition is associative.  Moreover, for
every object $x$ there is a morphism $1_x \maps x \to x$ acting
as the identity for composition.  It follows that all the
structure of a category with a single object $x$ can be
summarized by saying that $\hom(x,x)$, is a {\it monoid}: a set
equipped with an associative product and identity element.
Loosely, we say that a category with only one object is a monoid.
Similarly, a category with only one object and all morphisms
invertible is a group.

Recall also that given two categories $C$ and $D$, a functor $F
\maps C \to D$ maps objects of $C$ to objects of $D$ and
morphisms of $C$ to morphisms of $D$ in a structure-preseving
manner.  In other words, $f \maps x \to y$ implies $F(f) \maps
F(x) \to F(y)$, and also $F(fg) = F(f)F(g)$ and $F(1_x) =
1_{F(x)}$.  Thinking of a group as a one-object category $G$, a
group representation is thus just a functor from $G$ to $\Vect$,
the category whose objects are (finite-dimensional) vector spaces
and morphisms are linear maps.   This leads us to define a
`representation' of a category $C$ to be a functor from $C$ to
$\Vect$.

An $n$-dimensional TQFT is a certain sort of representation of
the category $n\Cob$ of $n$-dimensional cobordisms.  This
category has compact oriented $(n-1)$-manifolds as objects, and
oriented cobordisms between such manifolds as morphisms.
Composition of cobordisms is given by gluing as shown in Figure
1.  The identity $1_M$ for any object $M$ is represented by the
cylinder $[0,1] \times M$.  In what follows we will often abuse
language and identify cobordisms with the manifolds with boundary
representing them, but it is important to keep in mind the
distinction: for example, composition is not strictly
associative, but only associative up to equivalence, unless we
treat cobordisms carefully \cite{SawinTQFT}.

\bfig
\[ 
 \xy % U-shaped cobordism
  (0,6)*{f};   %THIS IS THE LABEL FOR THIS COBORDISM, YOU CAN MOVE 
%BY CHANGING THE COORDINATES!
  (-3,1.5)*\ellipse(3,1){-};
  (3,1.5)*\ellipse(3,1){-};
  (-3,3)*{}="X1";
  (3,3)*{}="X2";
    "X1";"X2" **\crv{(-3,-2) & (3,-2)};
  (-9,3)*{}="X1";
  (9,3)*{}="X2";
    "X1";"X2" **\crv{(-9,-9) & (9,-9)};
 \endxy
\qquad  
\circ 
\qquad
%Fork shaped cobordism
 \xy
  (-6,5)*{g};   %THIS IS THE LABEL FOR THIS COBORDISM, 
%YOU CAN MOVE BY CHANGING THE COORDINATES!
  (0,9)*\xycircle(3,1){-};
  (-3,-3)*\ellipse(3,1){.};
  (3,-3)*\ellipse(3,1){.};
  (-3,-3)*\ellipse(3,1)__,=:a(-180){-};
  (3,-3)*\ellipse(3,1)__,=:a(-180){-};
  (-3,-6)*{}="1";
  (3,-6)*{}="2";
  (-9,-6)*{}="A2";
  (9,-6)*{}="B2";
    "1";"2" **\crv{(-3,-1) & (3,-1)};
  (-3,9)*{}="A";
  (3,9)*{}="B";
  (-3,3)*{}="A1";
  (3,3)*{}="B1";
   "A";"A1" **\dir{-};
   "B";"B1" **\dir{-};
    "B1";"B2" **\crv{(8,1)};
    "A1";"A2" **\crv{(-8,1)};
 \endxy
\qquad 
= 
\qquad
 %Gluing the fork and U
 \xy
  (-10,5)*{fg};   %THIS IS THE LABEL FOR THIS COBORDISM, 
%YOU CAN MOVE BY CHANGING THE COORDINATES!
  (0,9)*\xycircle(3,1){-};
  (-3,-3)*\ellipse(3,1){.};
  (3,-3)*\ellipse(3,1){.};
  (-3,-3)*\ellipse(3,1)__,=:a(-180){-};
  (3,-3)*\ellipse(3,1)__,=:a(-180){-};
  (-3,-6)*{}="1";
  (3,-6)*{}="2";
  (-9,-6)*{}="A2";
  (9,-6)*{}="B2";
    "1";"2" **\crv{(-3,-1) & (3,-1)};
  (-3,9)*{}="A";
  (3,9)*{}="B";
  (-3,3)*{}="A1";
  (3,3)*{}="B1";
   "A";"A1" **\dir{-};
   "B";"B1" **\dir{-};
    "B1";"B2" **\crv{(8,1)};
    "A1";"A2" **\crv{(-8,1)};
  (-3,-6)*{}="X1";
  (3,-6)*{}="X2";
    "X1";"X2" **\crv{(-3,-11) & (3,-11)};
  (-9,-6)*{}="X1";
  (9,-6)*{}="X2";
    "X1";"X2" **\crv{(-9,-18) & (9,-18)};
 \endxy
\]
1.  Composition in $n\Cob$
\efig

A representation of $n\Cob$ is thus a functor $Z\maps n\Cob \to
\Vect$.   The fact that $Z$ assigns to the cylindrical spacetime
$[0,1] \times M$ the identity on $Z(M)$, that is, the trivial
time evolution operator, corresponds to the Wheeler-DeWitt
equation.  The dynamics of the theory only becomes evident upon
considering nontrivial cobordisms.

The category $n\Cob$ is not merely a category; it has extra
structures in common with $\Vect$, and the definition of a TQFT
requires that the functor $Z \maps n\Cob \to \Vect$ preserve
these extra structures.   In fact, these extra structures are
important clues about the nature of higher-dimensional algebra.

First, both $n\Cob$ and $\Vect$ are `monoidal' categories.  For precise
definitions of this and other terms from category theory, see Mac Lane
\cite{maclane}; roughly speaking, a category is monoidal if it has
tensor products of objects and morphisms satisfying all the usual
axioms, and an object $1$ playing the role of identity for the tensor
product.  In $n\Cob$, the tensor product is given by disjoint union, as
shown in Figure 2, and the identity is the empty set.  In $\Vect$, the
tensor product is the usual tensor product of vector spaces, which has
$\C$ as its identity.  In a TQFT, the functor $Z \maps n\Cob \to \Vect$
is required to be `monoidal', that is, to preserve tensor products and
to send the identity object in $n\Cob$ to the identity object in
$\Vect$.

\bfig
\[
 \xy  %Fork shaped cobordism
    (0,15)*{f};
    (0,9)*\xycircle(3,1){-};
  (-3,-3)*\ellipse(3,1){.};
  (3,-3)*\ellipse(3,1){.};
  (-3,-3)*\ellipse(3,1)__,=:a(-180){-};
  (3,-3)*\ellipse(3,1)__,=:a(-180){-};
  (-3,-6)*{}="1";
  (3,-6)*{}="2";
  (-9,-6)*{}="A2";
  (9,-6)*{}="B2";
    "1";"2" **\crv{(-3,-1) & (3,-1)};
  (-3,9)*{}="A";
  (3,9)*{}="B";
  (-3,3)*{}="A1";
  (3,3)*{}="B1";
   "A";"A1" **\dir{-};
   "B";"B1" **\dir{-};
    "B1";"B2" **\crv{(8,1)};
    "A1";"A2" **\crv{(-8,1)};
 \endxy
\quad
\tensor
\quad
 \xy
    (20,15)*{g};
  (7,4.5)*\ellipse(3,1){-};
  (13,4.5)*\ellipse(3,1){-};
  (17,9)*{}="X1";
  (23,9)*{}="X2";
    "X1";"X2" **\crv{(17,4) & (23,4)};
  (11,9)*{}="X1";
  (29,9)*{}="X2";
    "X1";"X2" **\crv{(11,-3) & (29,-3)};
\endxy
\qquad
=
\qquad
  \xy  %Fork shaped cobordism
    (9,15)*{f \tensor g};
    (0,9)*\xycircle(3,1){-};
  (-3,-3)*\ellipse(3,1){.};
  (3,-3)*\ellipse(3,1){.};
  (-3,-3)*\ellipse(3,1)__,=:a(-180){-};
  (3,-3)*\ellipse(3,1)__,=:a(-180){-};
  (-3,-6)*{}="1";
  (3,-6)*{}="2";
  (-9,-6)*{}="A2";
  (9,-6)*{}="B2";
    "1";"2" **\crv{(-3,-1) & (3,-1)};
  (-3,9)*{}="A";
  (3,9)*{}="B";
  (-3,3)*{}="A1";
  (3,3)*{}="B1";
   "A";"A1" **\dir{-};
   "B";"B1" **\dir{-};
    "B1";"B2" **\crv{(8,1)};
    "A1";"A2" **\crv{(-8,1)};
  (7,4.5)*\ellipse(3,1){-};
  (13,4.5)*\ellipse(3,1){-};
  (17,9)*{}="X1";
  (23,9)*{}="X2";
    "X1";"X2" **\crv{(17,4) & (23,4)};
  (11,9)*{}="X1";
  (29,9)*{}="X2";
    "X1";"X2" **\crv{(11,-3) & (29,-3)};
 \endxy
\]
2.  Tensor product in $n\Cob$
\efig

Second, both $n\Cob$ and $\Vect$ are `symmetric' monoidal
categories.  In a symmetric monoidal category, there is for
any pair of objects $x,y$ a natural isomorphism, the `braiding',
\[          B_{x,y} \maps x \tensor y \to y \tensor x, \]
which is required to satisfy various axioms including the symmetry
equation
\[        B_{y,x} B_{x,y}  = 1_{x \tensor y} .\]
In $n\Cob$, the symmetry $B_{x,y}$ is a cobordism of the sort
shown in Figure 3. In $\Vect$, the symmetry is the usual
isomorphism of vector spaces $x \tensor y$ and $y \tensor x$.

\bfig
\[
B_{x,y} 
\qquad 
= 
\qquad
 \xy 
  (4,4)*\ellipse(3,1){-};
  (-4,4)*\ellipse(3,1){-};
  (4,-4)*\ellipse(3,1){.};
  (-4,-4)*\ellipse(3,1){.};
  (-4,-4)*\ellipse(3,1)__,=:a(-180){-};
  (4,-4)*\ellipse(3,1)__,=:a(-180){-};
  (-8,11)*{x};
  (8,11)*{y};
  (-8,-11)*{y};
  (8,-11)*{x};
  (-12,0)*{}; %label for cobordism can go here
  (11,8)*{}="TLL";
  (5,8)*{}="TL";
  (11,-8)*{}="BLL";
  (5,-8)*{}="BL";
  (-10.75,7.5)*{}="TRR";
  (-5,8)*{}="TR";
  (-11,-8)*{}="BRR";
  (-5,-8)*{}="BR";
     (0,3)*{}="T";
   (-3,0)*{}="R";
   (3,0)*{}="L";
   (0,-3)*{}="B";
   "TRR";"R" **\dir{-};
   "BL";"B" **\dir{-};
   "TR";"T" **\dir{-};
   "BLL";"L" **\dir{-};
   "TLL";"BR" **\dir{-};
   "TL";"BRR" **\dir{-};
 \endxy
\]
3. Symmetry in $n\Cob$
\efig

Third, both $n\Cob$ and $\Vect$ are `rigid' monoidal
categories.   These are
monoidal categories in which every object $x$ has a `dual'
$x^\ast$, and there are `unit' and `counit' maps
\be          i_x \maps 1 \to x \tensor x^\ast, \quad
             e_x \maps x \tensor x^\ast \to 1  \label{unitcounit} \ee
satisfying various axioms including the triangle identities, which say
that the following diagrams must commute:
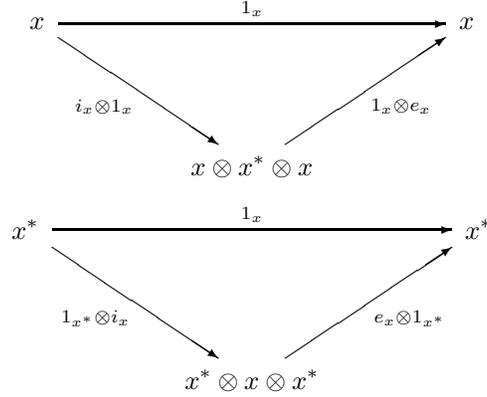
\[
\begin{diagram} [x \tensor x^\ast \tensor x]
\node{x} \arrow[2]{e,t}{1_x} \arrow{se,b}{i_x \tensor 1_x}
\node[2]{x}   \\
\node[2]{x \tensor x^\ast \tensor x} \arrow{ne,r}{1_x \tensor e_x}
\end{diagram}
\]
\[
\begin{diagram} [x^\ast \tensor x \tensor x^\ast]
\node{x^\ast} \arrow[2]{e,t}{1_x} \arrow{se,b}{1_{x^\ast} \tensor i_x}
\node[2]{x^\ast}   \\
\node[2]{x^\ast \tensor x \tensor x^\ast} \arrow{ne,r}{e_x \tensor 1_{x^\ast}}
\end{diagram}
\]
(For ease of exposition, we demand one extra axiom besides the usual ones,
namely that in the symmetric case the natural
morphism from $x$ to $x^{\ast\ast}$ corresponding to the `twist'
in Figure 31 be the identity.)
In $n\Cob$, $x^\ast$ is the manifold $x$ equipped with the
opposite orientation.  The unit and counit are the cylinders
shown in Figure 4, where we also depict the triangle identities.

\bfig
\[
e_x 
\quad 
=
\quad
 \xy 0;/r.16pc/: % U-shaped cobordism
  (-6,6.5)*{x^{\ast}};
  (6,6.5)*{x};
  (-11,0)*{}; %label for cobordism can go here
  (-3,1.5)*\ellipse(3,1){-};
  (3,1.5)*\ellipse(3,1){-};
  (-3,3)*{}="X1";
  (3,3)*{}="X2";
    "X1";"X2" **\crv{(-3,-2) & (3,-2)};
  (-9,3)*{}="X1";
  (9,3)*{}="X2";
    "X1";"X2" **\crv{(-9,-9) & (9,-9)};
 \endxy
\qquad \qquad
i_x 
\quad
= 
\quad
  \xy 0;/r.16pc/: % Upsidedown U-shaped cobordism
  (-6,-6.5)*{x};
  (6,-6.5)*{x^{\ast}};
  (-11,0)*{}; %label for cobordism can go here
  (-3,-1.5)*\ellipse(3,1){.};
  (3,-1.5)*\ellipse(3,1){.};
  (-3,-1.5)*\ellipse(3,1)__,=:a(-180){-};
  (3,-1.5)*\ellipse(3,1)__,=:a(-180){-};
  (-3,-3)*{}="X1";
  (3,-3)*{}="X2";
    "X1";"X2" **\crv{(-3,2) & (3,2)};
  (-9,-3)*{}="X1";
  (9,-3)*{}="X2";
    "X1";"X2" **\crv{(-9,9) & (9,9)};
 \endxy
\]
\vskip 1em
\[
 \xy 0;/r.14pc/: % zig-zagged tube with seams from gluing I
  (-3,-1.5)*\ellipse(3,1){.};
  (3,-1.5)*\ellipse(3,1){.};
  (-3,-1.5)*\ellipse(3,1)__,=:a(-180){-};
  (-3,-8)*\ellipse(3,1)__,=:a(-180){-};
  (3,-1.5)*\ellipse(3,1)__,=:a(-180){-};
  (9,-1.5)*\ellipse(3,1)__,=:a(-180){-};
  (-3,-3)*{}="X1";
  (3,-3)*{}="X2";
    "X1";"X2" **\crv{(-3,2) & (3,2)};
  (-9,-3)*{}="Y1";
  (9,-3)*{}="Y2";
    "Y1";"Y2" **\crv{(-9,9) & (9,9)};
  (9,-1.5)*\ellipse(3,1){.};
  (9,-3)*{}="Z1";
  (15,-3)*{}="Z2";
    "Z1";"Z2" **\crv{(9,-8) & (15,-8)};
  (3,-3)*{}="W1";
  (21,-3)*{}="W2";
    "W1";"W2" **\crv{(3,-15) & (21,-15)};
  (9,6)*\ellipse(3,1){-};
  (-3,-8)*\ellipse(3,1){.};
  (15,12)*{}="TL";
  (21,12)*{}="TR";
  (-9,-16)*{}="BL";
  (-3,-16)*{}="BR";
  "TL";"Z2" **\dir{-};
  "TR";"W2" **\dir{-};
  "BL";"Y1" **\dir{-};
  "BR";"X1" **\dir{-};
  (18,15.5)*{x};
  (-6,-19.5)*{x};
 \endxy
\quad = \quad
 \xy 0;/r.14pc/: %Straight tube cobordism
  (0,6)*\ellipse(3,1){-};
  (0,-8)*\ellipse(3,1){.};
  (0,-8)*\ellipse(3,1)__,=:a(-180){-};
  (-3,12)*{}="TL";
  (3,12)*{}="TR";
  (-3,-16)*{}="BL";
  (3,-16)*{}="BR";
   "TL"; "BL" **\dir{-};
   "TR"; "BR" **\dir{-};
   (,15.5)*{x};
  (,-19.5)*{x};
\endxy
\qquad \quad
 \xy 0;/r.14pc/: % zig-zaged tube with seams from gluing II
  (3,-1.5)*\ellipse(3,1){.};
  (-3,-1.5)*\ellipse(3,1){.};
  (3,-1.5)*\ellipse(3,1)__,=:a(-180){-};
  (3,-8)*\ellipse(3,1)__,=:a(-180){-};
  (-3,-1.5)*\ellipse(3,1)__,=:a(-180){-};
  (-9,-1.5)*\ellipse(3,1)__,=:a(-180){-};
  (3,-3)*{}="X1";
  (-3,-3)*{}="X2";
    "X1";"X2" **\crv{(3,2) & (-3,2)};
  (9,-3)*{}="Y1";
  (-9,-3)*{}="Y2";
    "Y1";"Y2" **\crv{(9,9) & (-9,9)};
  (-9,-1.5)*\ellipse(3,1){.};
  (-9,-3)*{}="Z1";
  (-15,-3)*{}="Z2";
    "Z1";"Z2" **\crv{(-9,-8) & (-15,-8)};
  (-3,-3)*{}="W1";
  (-21,-3)*{}="W2";
    "W1";"W2" **\crv{(-3,-15) & (-21,-15)};
  (-9,6)*\ellipse(3,1){-};
  (3,-8)*\ellipse(3,1){.};
  (-15,12)*{}="TL";
  (-21,12)*{}="TR";
  (9,-16)*{}="BL";
  (3,-16)*{}="BR";
  "TL";"Z2" **\dir{-};
  "TR";"W2" **\dir{-};
  "BL";"Y1" **\dir{-};
  "BR";"X1" **\dir{-};
  (-18,15.5)*{x^{\ast}};
  (6,-19.5)*{x^{\ast}};
 \endxy
\quad = \quad
 \xy 0;/r.14pc/:%Straight tube cobordism
  (0,6)*\ellipse(3,1){-};
  (0,-8)*\ellipse(3,1){.};
  (0,-8)*\ellipse(3,1)__,=:a(-180){-};
  (-3,12)*{}="TL";
  (3,12)*{}="TR";
  (-3,-16)*{}="BL";
  (3,-16)*{}="BR";
   "TL"; "BL" **\dir{-};
   "TR"; "BR" **\dir{-};
   (,15.5)*{x^{\ast}};
  (,-19.5)*{x^{\ast}};
\endxy
\]
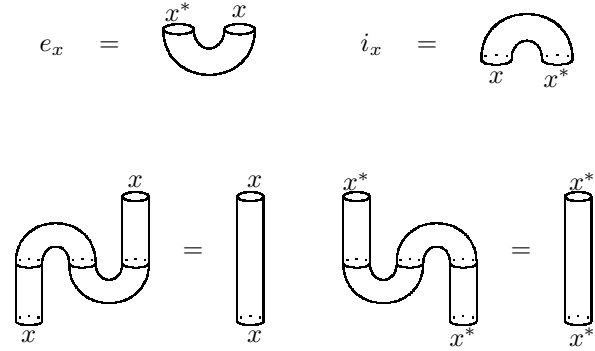
4.  Unit, counit, and triangle identities in $n\Cob$
\efig

In short, a TQFT is a rigid symmetric monoidal functor $X \maps
n\Cob \to \Vect$, that is, one preserving the rigid symmetric
monoidal structure. Now the category $\Hilb$, whose objects are
(finite-dimensional) Hilbert spaces and whose morphisms are linear
maps, is also rigid symmetric monoidal.  A `unitary' TQFT is a
rigid symmetric monoidal functor $Z\maps n\Cob \to \Hilb$ which
is also compatible with a second  sort of duality structure.
The operation $\dagger \maps n\Cob \to n\Cob$ taking each object
to itself and taking each cobordism $f \maps x \to y$ to the
orientation-reversed cobordism $f^\dagger \maps y \to x$ is a
contravariant functor, that is, $1_x^\dagger = 1_x$  and
$(fg)^\dagger = g^\dagger f^\dagger$.   There is also a
contravariant functor $\dagger \maps \Hilb \to \Hilb$ taking each
object to itself and taking each linear map $f\maps x \to y$ to
the Hilbert space adjoint $f^\dagger \maps y \to x$.   A unitary
TQFT  must satisfy $Z(f^\dagger) = Z(f)^\dagger$ for all
morphisms $f$.   Given cobordisms $f$ and $g$ from the empty set
to $x$, the inner product of the vectors $Z(f)1$ and $Z(g)1$ is
then given by $Z(f^\dagger g)1$.  If $Z(x)$ is spanned by vectors
of this form, the inner product in $Z(x)$ is thus determined by $Z$.

To conclude, it is crucial to note that Atiyah never propounded the
notion of a TQFT as a panacea for the problems of quantum field
theories without prior geometry.  Indeed, the TQFTs we understand
so far appear to reduce in the classical limit to field theories
with no local degrees of freedom, that is, for which all
solutions are {\it locally} physically equivalent.  One does not
expect a realistic theory of quantum gravity to have this
property.    Indeed, while there are many clues indicating that
quantum gravity is closely {\it related} to known TQFTs
\cite{BaezKQG}, it may be an inherently more complex sort
of theory.   Ideas from higher-dimensional algebra, however, may
still be very useful \cite{Baez}.

\section{Cobordisms}

While it is customary to begin in field theory by writing down a
Lagrangian, the most efficient ways to construct TQFTs tend to be
algebraic in flavor.   Since a TQFT is a rigid symmetric monoidal
functor from $n\Cob$ to $\Vect$, one can begin by
describing $n\Cob$ as a rigid symmetric monoidal category in
terms of generators and morphisms.  Here the `generators' are
morphisms from which one can obtain {\it all} the morphisms by
the operations present:  composition, tensor product, the
symmetry, and duals.  Then, to actually construct a TQFT, one
merely needs to assign objects and morphisms in $\Vect$ to all of
the generators of $n\Cob$, and check that the relations hold.

How does one determine generators and relations for $n\Cob$?
Assuming momentarily that we already understand the {\it objects}
in $n\Cob$ and their automorphisms, we can obtain generators
for the remaining {\it morphisms} using Morse
theory \cite{Milnor2}.  A cobordism from $M$
to $M'$ can be represented as an $n$-manifold $N$ having boundary
identified with $\overline M \cup M$.  If we put a `height'
function on $N$ --- a smooth real function $F$ with  $F|_M = 0$
and $F|_{M'} = 1$ --- generically it will be a Morse function.
That is, it will have only nondegenerate critical points $p_i$,
occuring at distinct levels $F(p_i) = t_i$.   Slicing $N$ along
level sets of $F$ between the critical levels $t_i$ amounts to
factoring our cobordism as a product of simple `generating'
cobordisms, as shown in Figure 5.

\bfig
\[
 \xy  
  (0,10)*\ellipse(3,1){-};
  (0,-12)*\ellipse(3,1)__,=:a(-180){-};
  (0,-12)*\ellipse(3,1){.};
  (0,7.5)*\ellipse(3,1){.};
  (0,7.5)*\ellipse(3,1)__,=:a(-180){-};
  (-3,0)*\ellipse(3,1){.};
  (3,0)*\ellipse(3,1){.};
  (-3,0)*\ellipse(3,1)__,=:a(-180){-};
  (3,0)*\ellipse(3,1)__,=:a(-180){-};
   (0,-6.5)*\ellipse(3,1){.};
  (0,-6.5)*\ellipse(3,1)__,=:a(-180){-};
  (-3,0)*{}="1";
  (3,0)*{}="2";
  (-9,0)*{}="xA2";
  (9,0)*{}="xB2";
    "1";"2" **\crv{(-3,5) & (3,5)};
  (-3,15)*{}="xA";
  (3,15)*{}="xB";
  (-3,9)*{}="xA1";
  (3,9)*{}="xB1";
   "xA";"xA1" **\dir{-};
   "xB";"xB1" **\dir{-};
    "xB1";"xB2" **\crv{(8,7)};
    "xA1";"xA2" **\crv{(-8,7)};
  (-3,0)*{}="1";
  (3,0)*{}="2";
  (-9,0)*{}="A2";
  (9,0)*{}="B2";
    "1";"2" **\crv{(-3,-5) & (3,-5)};
  (-3,-13)*{}="A";
  (3,-13)*{}="B";
  (-3,-9)*{}="A1";
  (3,-9)*{}="B1";
   "A";"A1" **\dir{-};
   "B";"B1" **\dir{-};
   "B1";"B2" **\crv{(8,-7)};
    "A1";"A2" **\crv{(-8,-7)};
    %TUBE TIP*********************
   (-3,20)*{}="f1";
   (-3,-19.5)*{}="f2";
   (3,20)*{}="f3";
   (3,-19.5)*{}="f4";
    (-3,-24)*{}="xf2";
   (3,-24)*{}="xf4";
   "f2";"xf2" **\dir{-};
   "f4";"xf4" **\dir{-};
   "f1";"xA" **\dir{-};
   "f3";"xB" **\dir{-};
   "f2";"A" **\dir{-}; %ADJUST HERE IF PLANES ARE NOT TRANSPERANT
   "f4";"B" **\dir{-};
   %End tube tips***************
    (0,-8)*{}="A"; %START OF SLICES
    (18,0)*{}="B";
    (0,8)*{}="G";
    (-18,0)*{}="Y";
    (-8,4.5)*{}="L";
    (8,4.5)*{}="R";
    "A";"B" **\dir{-};
    "Y";"L" **\dir{-};
    "Y";"A" **\dir{-};
    "R";"B" **\dir{-};
    %New SLICE
    (0,7)*{}="A"; %START OF SLICES
    (18,15)*{}="B";
    (0,23)*{}="G";
    (-18,15)*{}="Y";
    "A";"B" **\dir{-};
    "Y";"G" **\dir{-};
    "Y";"A" **\dir{-};
    "B";"G" **\dir{-};
        (0,-21)*{}="A"; %START OF SLICES
    (18,-13)*{}="B";
    (0,-5)*{}="G";
    (-18,-13)*{}="Y";
    (6,-8)*{}="R";
     (-6,-8)*{}="L";
    "A";"B" **\dir{-};
    "Y";"L" **\dir{-};
    "Y";"A" **\dir{-};
    "R";"B" **\dir{-};
 \endxy
\qquad \qquad \qquad
 \xy
 (0,7.5)*\ellipse(3,1){-};
  (-3,0)*\ellipse(3,1){-};
  (3,0)*\ellipse(3,1){-};
   (0,-6.5)*\ellipse(3,1){-};
   (0,9.5)*{}="1";
   (0,5.5)*{}="2";
   (0,-5)*{}="3";
   (0,-9)*{}="4";
   {\ar@{=>} "1";"2"};
   {\ar@{=>} "3";"4"};
 \endxy
\]
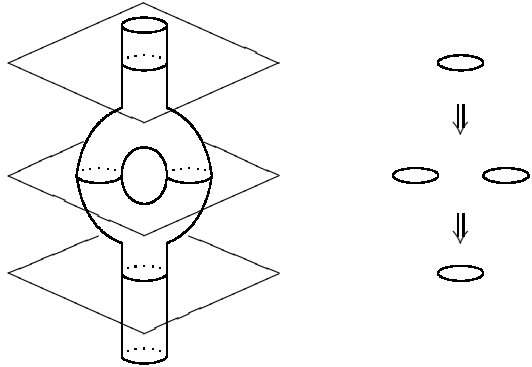
5. Describing a cobordism using Morse theory
\efig
\noindent
We can visualize the result as a `movie' of $N$ in which each
`frame' is obtained from the previous one by attaching a
$j$-handle --- that is, cutting out a copy of $D^{n-j-1} \times
S^{j-1}$ and gluing in a copy of $S^{n-j} \times D^{j}$ --- where
$j$ is the number of negative eigenvalues of the Hessian of  $F$
at the intervening critical point.    These basic processes are
shown for $n = 2$ in Figure 6; the cases $j = 0,1,2$ are called
the birth of a circle, the death of a circle, and the saddle
point, respectively.

\bfig
\[ 
 \vcenter{\xy
    (0,-4)*\ellipse(5,2){.};
    (0,-4)*\ellipse(5,2)__,=:a(-180){-};
    (-5,-8)*{}="TL";
    (5,-8)*{}="TR";
     "TL";"TR" **\crv{(-5,4) & (5,4) };
     (0,1)*{\bullet};
 \endxy}
\qquad \qquad
 \vcenter{\xy
    (0,4)*\ellipse(5,2){-};
    (0,18)*{};
    (-5,8)*{}="TL";
    (5,8)*{}="TR";
     "TL";"TR" **\crv{(-5,-4) & (5,-4) };
     (0,-1)*{\bullet};
 \endxy}
\qquad \qquad
 \vcenter{\xy
  (20,2)*{}="RU";
  (16,-3)*{}="RD";
  (-16,2)*{}="LU";
  (-20,-3)*{}="LD";
  "RU";"RD" **\crv{(4,2) & (4,-1)};
  "LD";"LU" **\crv{(-4,-2) & (-4,1)};
    (7.5,0)*{}="x1";
    (-7.5,0)*{}="x2";
     "x1"; "x2" **\crv{(7,-10) & (-7,-10)};
  (16,-20)*{}="RDD";
  (-20,-20)*{}="LDD";
   (20,-15)*{}="RUD";
   (-16,-15)*{}="LUD";
   (-16,-2.5)*{}="A";
   (16,-15)*{}="B";
        "RD"; "RDD" **\dir{-};
        "LD"; "LDD" **\dir{-};
        "A"; "LUD" **\dir{.};
        "RDD"; "LDD" **\dir{-};
        "RU"; "RUD" **\dir{-};
        "LU"; "A" **\dir{-};
        "B"; "RUD" **\dir{-};
        "B"; "LUD" **\dir{.};
        (0,-8)*{\bullet};
 \endxy }
\]
6.  Handle attachments for $n = 2$
\efig

Of course, the manifold $N$ admits many different Morse functions, so
the cobordism it represents can be expressed in many different ways as
the product of a series of handle attachments.  However, given two Morse
functions $F_0$ and $F_1$, we can interpolate between them by a smooth
family of functions $F_s$.  Generically, the functions $F_s$ will be
Morse functions except for finitely many values $s_i$ at which the level
of one critical point passes another, two critical points coalesce, or a
critical point splits in two.  The study of these generic paths between
Morse functions is known as Cerf theory \cite{Cerf,Kerler}.  In the same
sense as which handle attachments give {\it generators} for $n\Cob$,
these paths between Morse functions give {\it relations,} known as
handle slides and cancellations.  We can visualize these as `movie
moves' going between two different movies of the same cobordism.  An
example of a handle cancellation for the $n = 2$ case is shown in Figure
7.

\bfig
\[
 \xy 
  (0,4)*\ellipse(3,1){-};
  (0,-7)*\ellipse(3,1){.};
  (0,-7)*\ellipse(3,1)__,=:a(-180){-};
  (-3,8)*{}="TL";
  (3,8)*{}="TR";
  (-3,-14)*{}="BL";
  (3,-14)*{}="BR";
  (3,0)*{}="MT";
  (3,-8)*{}="MB";
  (7,0)*{}="2";
  (12,0)*{}="3";
   "TL"; "BL" **\dir{-};
   "TR"; "MT" **\dir{-};
   "BR"; "MB" **\dir{-};
   "MT";"2" **\crv{(3,-4)& (7,-4)};
   "2";"3" **\crv{(7,6)& (12,6)};
   "MB";"3" **\crv{(4,-4)& (12,-7)};
\endxy
\qquad = \qquad
 \xy %Straight tube cobordism
  (0,4)*\ellipse(3,1){-};
  (0,-7)*\ellipse(3,1){.};
  (0,-7)*\ellipse(3,1)__,=:a(-180){-};
  (-3,8)*{}="TL";
  (3,8)*{}="TR";
  (-3,-14)*{}="BL";
  (3,-14)*{}="BR";
   "TL"; "BL" **\dir{-};
   "TR"; "BR" **\dir{-};
\endxy
\]
7.  A handle cancellation for $n = 2$
\efig

For the case $n =1$ it is easy to use these ideas to give
a purely algebraic description of $n\Cob$.  It is (up to
the standard notion of equivalence of categories) just the free
rigid symmetric monoidal category on one object $x$!   The
object $x$ corresponds to the positively oriented point.  As
shown in Figure 8, the unit and counit
\[   i_x \maps 1 \to x \tensor x^\ast, \qquad
     e_x \maps x^\ast \tensor x \to 1 ,\]
correspond to the
two types of handle attachments, namely the birth and the death of an $S^0$
(a pair of oppositely oriented points).

\bfig
\[ 
i_x 
\quad 
= 
\quad
\vcenter{
 \xy
(-4,-3.5)*{x}; (4,-3.5)*{x^{\ast}}; (-8,3)*{};  % can put label here
 (-4,0)*{\bullet}; (4,0)*{\bullet} **\crv{(-4.5,12) & (4.5,12)};
 \endxy}
\qquad 
\qquad
e_x 
\quad
=
\quad 
\vcenter{
 \xy
(4,3.5)*{x}; (-4,3.5)*{x^{\ast}}; (-8,-3)*{};  % can put label here
 (4,0)*{\bullet}; (-4,0)*{\bullet} **\crv{(4.5,-12) & (-4.5,-12)};
 \endxy}
\]
8.  Handle attachments for $n = 1$
\efig

\noindent Similarly, as shown in Figure 9, the triangle
identities correspond to handle cancellations.

\bfig
\[
 \xy 
  (-6,10)*{\bullet}="1E";
  (-6,0)*{}="1";
  (0,0)*{}="2";
  (6,0)*{}="3";
  (6,-10)*{\bullet}="3B";
    "2";"1" **\crv{(0,-6)& (-6,-6)};
    "3";"2" **\crv{(6,6)& (0,6)};
    "1";"1E" **\dir{-};
    "3B";"3" **\dir{-};
    (-6,13.5)*{x^{\ast}};
    (6,-13.5)*{x^{\ast}};
 \endxy
 \quad = \quad
 \xy
 (0,10)*{\bullet};
 (0,-10)*{\bullet};
 **\dir{-};
 (0,13.5)*{x^{\ast}};
 (0,-13.5)*{x^{\ast}};
 \endxy
\qquad \qquad \qquad
 \xy
  (6,10)*{\bullet}="1E";
  (6,0)*{}="1";
  (0,0)*{}="2";
  (-6,0)*{}="3";
  (-6,-10)*{\bullet}="3B";
    "2";"1" **\crv{(0,-6)& (6,-6)};
    "3";"2" **\crv{(-6,6)& (0,6)};
    "1";"1E" **\dir{-};
    "3B";"3" **\dir{-};
    (6,13.5)*{x};
    (-6,-13.5)*{x};
 \endxy
 \quad = \quad
 \xy
 (0,10)*{\bullet};
 (0,-10)*{\bullet};
 **\dir{-};
  (0,13.5)*{x};
 (0,-13.5)*{x};
 \endxy
\]
9.  Handle cancellations for $n = 1$
\efig

This simple result indicates that the rigid symmetric monoidal
structure captures all the essential aspects of $n\Cob$ for $n =1$.  It
is tempting to seek similar purely algebraic presentations of $n\Cob$ in
higher dimensions.  For $n = 2$ we can achieve this using the principle
of `internalization'.  In its simplest form, this amounts to the fact
that any algebraic structure definable using commutative diagrams in the
category $\Set$ can be generalized to categories sufficiently resembling
$\Set$.  For example, a monoid can be defined as an object $x$ in $\Set$
equipped with a product $m \maps x \times x \to x$ and unit $i \maps 1
\to x$ making various diagrams commute.  Here $1$ denotes any
one-element set, and we use the standard trick of thinking of the
identity element of $x$ as the image of a map $i \maps 1 \to x$.  We can
generalize the definition to any monoidal category $C$, replacing
$\times$ with the tensor product in $C$ and $1$ with the identity object
of $C$, thus obtaining the notion of a `monoid object' in $C$.  For
example, a monoid object in $\Vect$ is an algebra.  It turns out that
$2\Cob$ is the `free rigid symmetric monoidal category on one
commutative monoid object with nondegenerate trace'.  The object in
question is $S^1$, and the product, identity, and trace $\tr \maps S^1
\to 1$ are shown in Figure 10, This result yields a complete
classification of 2-dimensional TQFTs \cite{DJ,Sawin2d}.

\bfig
\[  m 
\quad 
=
\quad
  \xy 0;/r.20pc/: 
    (0,-9)*\xycircle(3,1){.};
    (0,-4.5)*\ellipse(3,1)__,=:a(180){-};
  (-3,3)*\ellipse(3,1){-};
  (3,3)*\ellipse(3,1){-};
  (-3,6)*{}="1";
  (3,6)*{}="2";
  (-9,6)*{}="A2";
  (9,6)*{}="B2";
    "1";"2" **\crv{(-3,1) & (3,1)};
  (-3,-9)*{}="A";
  (3,-9)*{}="B";
  (-3,-3)*{}="A1";
  (3,-3)*{}="B1";
   "A";"A1" **\dir{-};
   "B";"B1" **\dir{-};
    "B1";"B2" **\crv{(8,1)};
    "A1";"A2" **\crv{(-8,1)};
    (-11.5,2)*{};           %CAN PUT LABEL HERE
 \endxy
\qquad \qquad
i 
\quad
= 
\quad
 \vcenter{\xy
    (0,0)*\ellipse(5,2){.};
    (0,0)*\ellipse(5,2)__,=:a(-180){-};
    (-5,0)*{}="TL";
    (5,0)*{}="TR";
     "TL";"TR" **\crv{(-5,12) & (5,12) };
     (-8,8)*{};               %CAN PUT LABEL HERE
 \endxy}
\qquad \qquad
\tr
\quad 
=
\quad
 \vcenter{\xy
    (0,4)*\ellipse(5,2){-};
    (0,18)*{};
    (-5,8)*{}="TL";
    (5,8)*{}="TR";
     "TL";"TR" **\crv{(-5,-4) & (5,-4) };
     (-8,8)*{};                %CAN PUT LABEL HERE
 \endxy}
 \]
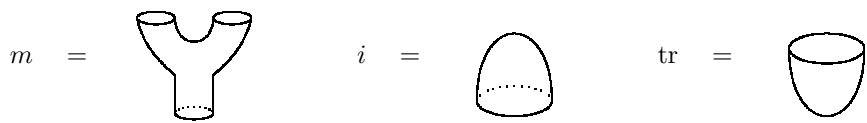
10.  $S^1$ as a commutative monoid object with nondegenerate
trace
\efig

Moving to higher dimensions, the best presentations of $3\Cob$ for the
purposes of constructing TQFTs are based on the Kirby calculus
\cite{Kerler,RT,SawinTQFT}.  While very algebraic in flavor, these
have not yet been distilled to a statement comparable to
those for $1\Cob$ and $2\Cob$.  The Kirby calculus also gives a
description of $4\Cob$ which has yielded a few TQFTs so far
\cite{Broda,Kirby}.  For $n \ge 6$ the theory of cobordisms becomes more
closely tied to homotopy theory, due to the $h$-cobordism theorem
\cite{Milnor}.  Also, while we will not go into it here, it is important
to note the existence of a theory of piecewise-linear (PL) manifolds
paralleling the smooth theory \cite{BW,CFS,CKY,FHK,Pachner}.  The
smooth and PL versions of $n\Cob$ are equivalent for $n \le 6$, but
not in general for larger $n$.

What we seek, however, is a unified algebraic framework for this
entire collection of results, one that applies to all dimensions
and explains the fascinating relationships between results in
neighboring dimensions.  Such a framework should clarify the
existing TQFT constructions, which appear at first to rely on
miraculous analogies between topology and algebra, and it should
aim at a classification of TQFTs.   We suggest that
$n$-category theory will provide such a framework.  The reason is
quite simple.  We sketched how to go about a `generators and
relations' description of the morphisms for $n\Cob$ assuming
we already had a description of the objects and their automorphisms, but
for higher $n$ these grow more complicated to describe as well.
Note however that in describing an $n$-dimensional manifold as a
`movie', each `frame', being an $(n-1)$-manifold, can itself be
regarded a `movie' of one lower dimension, as shown in Figure 11.
In other words, each object in $n\Cob$ gives a morphism in
$(n-1)\Cob$.   A fully algebraic description of $n\Cob$ should
therefore involve objects, morphisms between objects, 2-morphisms
between morphisms, and so on, up to $n$-morphisms.  We may
loosely call any such structure an $n$-category.

\bfig
\[
 \xy 
  (0,-12)*\ellipse(3,1)__,=:a(-180){-};
  (0,-12)*\ellipse(3,1){.};
  (-3,0)*\ellipse(3,1){.};
  (3,0)*\ellipse(3,1){.};
  (-3,0)*\ellipse(3,1)__,=:a(-180){-};
  (3,0)*\ellipse(3,1)__,=:a(-180){-};
    (-3,4)*\ellipse(3,1){-};
  (3,4)*\ellipse(3,1){-};
   (0,-6.5)*\ellipse(3,1){.};
  (0,-6.5)*\ellipse(3,1)__,=:a(-180){-};
  (-3,0)*{}="1";
  (3,0)*{}="2";
  (-9,0)*{}="A2";
  (9,0)*{}="B2";
    "1";"2" **\crv{(-3,-5) & (3,-5)};
  (-3,-13)*{}="A";
  (3,-13)*{}="B";
  (-3,-9)*{}="A1";
  (3,-9)*{}="B1";
   "A";"A1" **\dir{-};
   "B";"B1" **\dir{-};
   "B1";"B2" **\crv{(8,-7)};
    "A1";"A2" **\crv{(-8,-7)};
    %TUBE TIP*********************
      (-9,8)*{}="tl";
      (-3,8)*{}="tr";
      (9,8)*{}="er";
      (3,8)*{}="el";
   (-3,-19.5)*{}="f2";
   (3,-19.5)*{}="f4";
    (-3,-24)*{}="xf2";
   (3,-24)*{}="xf4";
   "1";"tr" **\dir{-};
   "2";"el" **\dir{-};
   "A2";"tl" **\dir{-};
   "B2";"er" **\dir{-};
   "f2";"xf2" **\dir{-};
   "f4";"xf4" **\dir{-};
   "f2";"A" **\dir{-}; %ADJUST HERE IF PLANES ARE NOT TRANSPERANT
   "f4";"B" **\dir{-};
   %End tube tips***************
    (0,-8)*{}="A"; %START OF SLICES
    (18,0)*{}="B";
    (0,8)*{}="G";
    (-18,0)*{}="Y";
    (-9.2,4)*{}="L";
    (9.2,4)*{}="R";
    (-3,6.5)*{}="LL";
    (3,6.5)*{}="RR";
    "A";"B" **\dir{-};
    "LL";"G" **\dir{-};
    "Y";"L" **\dir{-};
    "Y";"A" **\dir{-};
    "G";"RR" **\dir{-};
    "R";"B" **\dir{-};
        (0,-21)*{}="A"; %START OF SLICES
    (18,-13)*{}="B";
    (0,-5)*{}="G";
    (-18,-13)*{}="Y";
    (6,-8)*{}="R";
     (-6,-8)*{}="L";
    "A";"B" **\dir{-};
    "Y";"L" **\dir{-};
    "Y";"A" **\dir{-};
    "R";"B" **\dir{-};
 \endxy
\qquad \qquad
 \xy
  (-3,0)*\ellipse(3,1){-};
  (3,0)*\ellipse(3,1){-};
   (0,-6.5)*\ellipse(3,1){-};
   (0,-2.5)*{}="d1";
   (0,-6.5)*{}="d2";
   {\ar@{=>} "d1";"d2"};
   (-6,4)*{}="d1";
   (-6,-4)*{}="d2";
   "d1";"d2" **\dir{-};
   (6,4)*{}="d1";
   (6,-4)*{}="d2";
   "d1";"d2" **\dir{-};
   (0,-9)*{}="d1";
   (0,-17)*{}="d2";
   "d1";"d2" **\dir{-};
\endxy
\qquad \qquad \xy
  (-6,1.5)*{\bullet};
  (-6,-1.5)*{\bullet};
  (6,1.5)*{\bullet};
  (6,-1.5)*{\bullet};
  (0,-11.5)*{\bullet};
  (0,-14.5)*{\bullet};
  (0,-2.5)*{}="d1";
   (0,-6.5)*{}="d2";
   {\ar@{=>} "d1";"d2"};
   (-8,0)*{}="d2";
   (-12,0)*{}="d1";
   {\ar@{->} "d1";"d2"};
   (-1.5,0)*{}="d2";
   (-5,0)*{}="d1";
   {\ar@{->} "d1";"d2"};
   (8,0)*{}="d1";
   (12,0)*{}="d2";
   {\ar@{->} "d1";"d2"};
   (1.5,0)*{}="d1";
   (5,0)*{}="d2";
   {\ar@{->} "d1";"d2"};
   (-6,-13)*{}="d1";
   (-2,-13)*{}="d2";
   {\ar@{->} "d1";"d2"};
   (2,-13)*{}="d1";
   (6,-13)*{}="d2";
   {\ar@{->} "d1";"d2"};
\endxy
\]
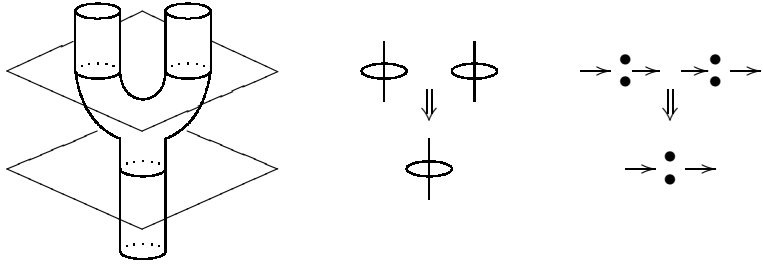
11.  Frames as movies of one lower dimension
\efig

If we take the objects of our $n$-category to be
0-manifolds, and the morphisms to be 1-manifolds with boundary,
the most general sort of 2-morphism will be some kind of
2-manifold with corners, as shown in Figure 12.

\bfig
\[
 \xy 
  (20,2)*{\bullet}="RU";
  (16,-3)*{\bullet}="RD";
  (-16,2)*{\bullet}="LU";
  (-20,-3)*{\bullet}="LD";
  "RU";"RD" **\crv{(4,2) & (4,-1)}; ?(.1)*\dir{<}; ?(.85)*\dir{<};
  "LD";"LU" **\crv{(-4,-2) & (-4,1)}; ?(.08)*\dir{<}; ?(.85)*\dir{<};
    (7.5,0)*{}="x1";
    (-7.5,0)*{}="x2";
     "x1"; "x2" **\crv{(7,-10) & (-7,-10)};
  (16,-20)*{\bullet}="RDD";
  (-20,-20)*{\bullet}="LDD";
   (20,-15)*{\bullet}="RUD";
   (-16,-15)*{\bullet}="LUD";
   (-16,-2.5)*{}="A";
   (16,-15)*{}="B";
        "RD"; "RDD" **\dir{-};
        "LD"; "LDD" **\dir{-};
        "A"; "LUD" **\dir{.};
        "RDD"; "LDD" **\dir{-}; ?(.55)*\dir{>};
        "RU"; "RUD" **\dir{-};
        "LU"; "A" **\dir{-};
        "B"; "RUD" **\dir{-};
        "B"; "LUD" **\dir{.}; ?(.45)*\dir{<};

 \endxy
\]
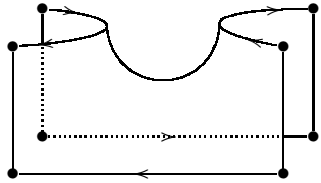
12. 2-manifold with corners as 2-morphism between 1-manifolds with boundary
\efig

\noindent
This is a considerable nuisance, since as $n$ increases it
becomes more and more difficult to specify a precise class of
`$n$-manifolds with corners' and precise recipes for composing
them.  However, this added complexity is of vital importance,
since it permits the definition of an `extended' TQFT, one which
behaves well under the extra cutting and pasting constructions
available in this context.   Heuristic reasoning involving path
integrals suggests that the TQFTs described in terms of local
Lagrangians should be of this extended sort, and so far this has
been borne out in rigorous work on important examples
\cite{Freed,FQ,Lawrence,Q,Walker}.  It is, in fact, the theory of
extended TQFTs that provides the best information about the
relationship between higher-dimensional algebra and TQFTs.

We could at this point attempt to define `manifolds with corners' more
precisely, and define composition operations on them.  Various
approaches have already been successfully pursued in the work
mentioned above.  However, in order to get an idea of what a convenient
formalism should eventually look like, we prefer to turn to the theory of
$n$-categories, and see what that suggests.

\section{Strict $n$-Categories}

Often when people refer to $n$-categories they mean what we shall call
`strict' $n$-categories.  These appear not to be sufficiently general
for TQFT applications, but unlike the more general `weak'
$n$-categories, they have already been defined for all $n$.  In what
follows we briefly sketch this definition and some of its implications,
while in the next section we indicate the importance of weakening it in
certain ways.  Readers familiar with strict $n$-categories can skip this
section.  In the rest of this section we omit the qualifier `strict'.

The most elegant approach involves the theory of `enriched'
categories \cite{Kelly}.  This is based on the observation that
in the definition of a category $C$, the category $\set$, whose
objects are sets and whose morphisms are functions, plays a
distinguished role.  The reason is that:
\vskip 1em
\noindent For every pair of objects $x,y$ in $C$ there is a {\it set}
$\,hom(x,y)$ of morphisms, and for every triple of objects $x,y,z$
in $C$ composition is a {\it function} $\circ \maps \hom(x,y)
\times \hom(y,z) \to \hom(x,z)$.
\vskip 1em \noindent
Note also that the monoidal structure of $\Set$, the
Cartesian product $\times$, plays a role here.  All the rest of
the category axioms can be written out as commutative
diagrams in the category $\set$, and all of these diagrams make
sense in any monoidal category.  Thus one can relativize the
definition of a category by letting an arbitrary monoidal
category $K$ play the role that $\Set$ does here.  In other
words, a category $C$ `enriched over $K$', or `$K$-category',
is a collection of objects for which:
\vskip 1em
\noindent For every pair of objects $x,y$ in $C$ there is a {\it
object} $\,\hom(x,y)$ in $K$, and for every triple of objects
$x,y,z$ in $C$ there is a {\it morphism} $\circ
\maps \hom(x,y) \tensor \hom(y,z) \to \hom(x,z)$ in $K$.  \vskip 1em
\noindent One also demands that the usual axioms of a category
hold, translated into commutative diagrams in $K$.

A simple example is the category $\Vect$, which is enriched over
itself, or `closed' \cite{EK}.  That is, given vector spaces $x$
and $y$, the set $\hom(x,y)$ is actually a {\it vector space},
and composition $\circ \maps \hom(x,y) \tensor \hom(y,z) \to
\hom(x,z)$ is actually a {\it linear map}.   Another example
is the category of modules over a ring, which
is enriched over the category of abelian groups.

This notion of enriched category permits a wonderful recursive
definition of $n$-categories, as follows.  We say a category is
`small' if the collection of objects is a set.  The category of
all small categories is denoted $\Cat$ --- note that the `smallness'
condition prevents Russell-type paradoxes.  Now $\Cat$ is actually a
monoidal category, with the identity $1$ taken as any category with
one object $x$ and one morphism $1_x$, and with the tensor product
being the usual Cartesian product $\times$ of categories.  This
product is an obvious generalization of the Cartesian product of sets,
as shown in Figure 13.

\bfig
\[
 \xy 
    (-16,-2)*{\bullet}; % categories C and D.
    (-16,6)*{\bullet};
    (0,-2)*{\bullet};
    (0,6)*{\bullet};
    (8,-2)*{\bullet};
    (8,6)*{\bullet};
    (16,-2)*{\bullet};
    (16,6)*{\bullet};
    (0,15)*{\bullet};
    (8,15)*{\bullet};
    (16,15)*{\bullet};
    (9,18)*{S};
    (-20,2.5)*{T};
    (25,2.5)*{S \times T};
 \endxy
\qquad \qquad
 \xy
    (-13,-2)*+{\bullet}="1";
    (-13,6)*+{\bullet}="2";
    (0,-2)*+{\bullet}="d";
    (0,6)*+{\bullet}="a";
    (10,-2)*+{\bullet}="e";
    (10,6)*+{\bullet}="b";
    (20,-2)*+{\bullet}="f";
    (20,6)*+{\bullet}="c";
    (0,15)*+{\bullet}="4";
    (10,15)*+{\bullet}="5";
    (20,15)*+{\bullet}="6";
    (10,18)*{C};
    (-18,2.5)*{D};
    (29,2.5)*{C \times D};
    {\ar@/^.25pc/ "2";"1"};
    {\ar@/^.25pc/ "a";"b"};
    {\ar@/^.25pc/ "b";"c"};
    {\ar@/^.25pc/ "d";"e"};
    {\ar@/^.25pc/ "e";"f"};
    {\ar@/^.25pc/ "a";"d"};
    {\ar@/^.25pc/ "b";"e"};
    {\ar@/^.25pc/ "c";"f"};
    {\ar@/^.25pc/ "4";"5"};
    {\ar@/^.25pc/ "5";"6"};
 \endxy
\]
\vskip 1em
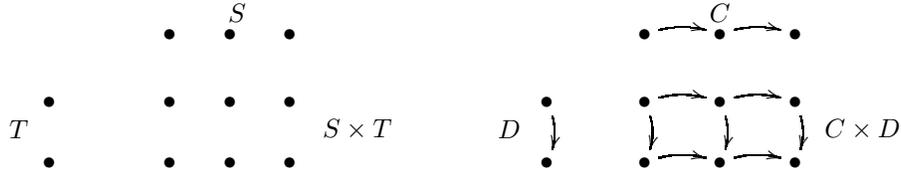
13.  The Cartesian product of sets $S$ and $T$, and of categories $C$ and $D$
\efig

\noindent  To be precise, the objects of $C \times D$, written
as $x \times y$, are just ordered pairs consisting of an object $x$ of
$C$ and an object $y$ of $D$.  The morphisms of $C \times D$
can be described using generators and relations.
Given a morphism $f \maps x \to x'$
in $C$ and a morphism $g \maps y \to y'$ in $D$, there
are morphisms
\[  x \times g \maps x \times y \to x \times y',\qquad
f \times y \maps x \times y \to x' \times y  \]
in $C \times D$.
These are the generators; the relations say that
\[    (f \times y)(f' \times y) = ff' \times y, \qquad
      (x \times g)(x \times g') = x \times gg', \]
and very importantly, that diagrams of the following form commute:
\be \label{cartprod}
\begin{diagram}
\node{x \times y} \arrow{e,t}{f \times y} \arrow{s,l}{x \times g}
\node{x' \times y} \arrow{s,r}{x' \times g}\\
\node{x \times y'}  \arrow{e,t}{f \times y'} \node{x' \times y'.}
\end{diagram}
\ee
This implies that all the squares in Figure 13 commute.

The definition of (strict) $n$-categories is then as follows.
2-categories are simply categories enriched over $\Cat$.  The
category $2\Cat$ of small 2-categories, in turn, has a Cartesian
product making it into a monoidal category.  This allows
us to define 3-categories as categories enriched over
$2\Cat$.  In general, $n\Cat$ is defined as the category of small
categories enriched over $(n-1)\Cat$, which is monoidal
when equipped with its Cartesian product.   The Cartesian product
at each stage is defined by a generalization of the Cartesian
product of categories to the enriched context \cite{Kelly}.

For the reader dizzied by the rapid ascent up this recursive
ladder, let us briefly pause to contemplate the case of
2-categories \cite{KS}.  Here for any pair of objects $x$ and
$y,$ $\hom(x,y)$ is a {\it category}.  Objects in $\hom(x,y)$ should be
thought of as morphisms from $x$ to $y$, while the morphisms in
$\hom(x,y)$ should be thought of as `morphisms between morphisms'
or 2-morphisms.   The 2-morphisms are sometimes drawn as
2-dimensional surfaces labelled with double arrows; in Figure 14
we show objects $x,y$, morphisms $f\maps x \to y$ and $g\maps x
\to y$, and a 2-morphism $\alpha \maps f \doublearrow g$.

\bfig
\[
 \xy 
  (-15,0)*{\bullet}="1"; %of a 1-morphism with a 2-morphism.
  (0,0)*{\bullet}="2";
 "1";"2" **\crv{(-12,9) & (-3,9)};
  "1";"2" **\crv{(-12,-9) & (-3,-9)};
    (-7.5,6.75)*{\scriptstyle >}+(0,3)*{\scriptstyle f};
    (-7.5,-6.75)*{\scriptstyle >}+(0,-3)*{\scriptstyle g};
    (-17.5,0)*{\scriptstyle x};
  (2.5,0)*{\scriptstyle y};
  {\ar@{=>}_{\scriptstyle \alpha}(-7.5,3)*{};(-7.5,-3)*{}} ;
 \endxy
\]
14.  Diagram of a 2-morphism $\alpha \maps f \doublearrow g$
\efig

Given objects $x,y$ in a 2-category $C$ and morphisms $f,g,h
\maps x \to y$, we can compose 2-morphisms $\alpha \maps g
\doublearrow h$ and $\beta \maps f \doublearrow g$ to obtain a
2-morphism $\alpha\beta \maps f \doublearrow h$.  This operation,
which is really just composition of morphisms in the category
$\hom(x,y)$, is often called `vertical' composition, for reasons
made clear by Figure 15.

\bfig
\[
 \xy 
  (-15,0)*{\bullet}="1";
  (0,0)*{\bullet}="2";
 "1";"2" **\crv{(-12,9) & (-3,9)};
 "1";"2" **\crv{};
  "1";"2" **\crv{(-12,-9) & (-3,-9)};
    (-7.5,6.75)*{\scriptstyle >}+(0,3)*{\scriptstyle f};
    (-7.5,-6.75)*{\scriptstyle >}+(0,-3)*{\scriptstyle h};
    (-17.5,0)*{\scriptstyle x};
  (2.5,0)*{\scriptstyle y};
    (-7.5,0)*{\scriptstyle >}+(2.5,1.5)*{\scriptstyle g};
  {\ar@{=>}_{\scriptstyle \beta}(-7.5,5)*{};(-7.5,1.25)*{}} ;
  {\ar@{=>}_{\scriptstyle \alpha}(-7.5,-1.25)*{};(-7.5,-5)*{}} ;
 \endxy
\]
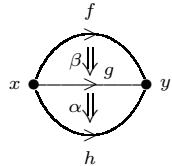
15.  Vertical composition of 2-morphisms
\efig

\noindent On the other hand, given objects $x,y$ and $z$, the
composition functor from $\hom(x,y) \times \hom(y,z)$ to
$\hom(x,z)$, gives various other operations.  Referring back to
the definition of the Cartesian product of categories,
one sees that composition takes an object in $\hom(x,y)$ and one
in $\hom(y,z)$ to one in $\hom(x,z)$: this is how one composes
1-morphisms in $C$.  Composition also takes an object in one of
these categories and a morphism in the other to a morphism in
$\hom(x,z)$.  This gives two ways to compose a 1-morphism and a
2-morphism in $C$ to obtain a 2-morphism, as shown in Figure 16.

\bfig
\[
 \xy 
  (-15,0)*{\bullet}="1";
  (0,0)*{\bullet}="2";
  (15,0)*{\bullet}="3";
 "2";"3" **\crv{(3,9) & (12,9)};
  "2";"3" **\crv{(3,-9) & (12,-9)};
  "1";"2" **\crv{};
    (7.5,6.75)*{\scriptstyle >}+(0,3)*{\scriptstyle g};
    (7.5,-6.75)*{\scriptstyle >}-(0,3)*{\scriptstyle g'};
    (-7.5,0)*{\scriptstyle >}+(0,3)*{\scriptstyle f};
  {\ar@{=>}_{\scriptstyle \beta}(7.5,3)*{};(7.5,-3)*{}} ;
    (-17.5,0)*{\scriptstyle x};
  (2.5,0)*{\scriptstyle y};
  (17.5,0)*{\scriptstyle z};
 \endxy
\qquad 
\qquad
 \xy 
  (-15,0)*{\bullet}="1";
  (0,0)*{\bullet}="2";
  (15,0)*{\bullet}="3";
 "1";"2" **\crv{(-12,9) & (-3,9)};
  "1";"2" **\crv{(-12,-9) & (-3,-9)};
  "2";"3" **\crv{};
    (-7.5,6.75)*{\scriptstyle >}+(0,3)*{\scriptstyle f};
    (-7.5,-6.75)*{\scriptstyle >}-(0,3)*{\scriptstyle f'};
    (7.5,0)*{\scriptstyle >}+(0,3)*{\scriptstyle g};
  {\ar@{=>}_{\scriptstyle \alpha}(-7.5,3)*{};(-7.5,-3)*{}} ;
    (-17.5,0)*{\scriptstyle x};
  (-2.5,0)*{\scriptstyle y};
  (17.5,0)*{\scriptstyle z};
 \endxy
\]
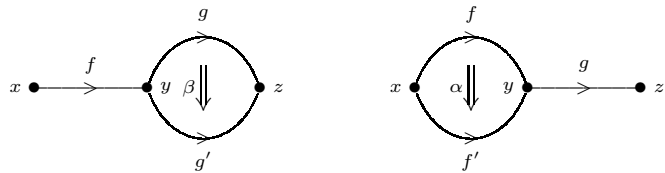
16.  Composition of a 1-morphism and a 2-morphism
\efig

\noindent  Thanks to eq.\ (\ref{cartprod}), we can
use these basic composition operations to define an operation
called `horizontal composition' of 2-morphisms, as shown in
Figure 17.  Given $f,f' \maps x \to y$, $g,g' \maps y \to z$,
$\alpha \maps f \doublearrow f'$ and $\beta \maps g \doublearrow
g'$, the horizontal composite of $\alpha$ and $\beta$ is a
2-morphism $\alpha \tensor \beta \maps gf \doublearrow g'f'$.

\bfig
\[
 \xy 
  (-15,0)*{\bullet}="1";
  (0,0)*{\bullet}="2";
  (15,0)*{\bullet}="3";
 "1";"2" **\crv{(-12,9) & (-3,9)};
  "1";"2" **\crv{(-12,-9) & (-3,-9)};
   "2";"3" **\crv{(3,9) & (12,9)} ;
    "2";"3" **\crv{ (3,-9)& (12,-9) };
    (-7.5,6.75)*{\scriptstyle >}+(0,3)*{\scriptstyle f};
    (7.5,6.75)*{\scriptstyle >}+(0,3)*{\scriptstyle g};
    (-7.5,-6.75)*{\scriptstyle >}+(0,-3)*{\scriptstyle f'};
    (7.5,-6.75)*{\scriptstyle >}+(0,-3)*{\scriptstyle g'};
    (-17.5,0)*{\scriptstyle x};
  (-2.5,0)*{\scriptstyle y};
  (17.5,0)*{\scriptstyle z};
  {\ar@{=>}_{\scriptstyle \beta}(7.5,3)*{};(7.5,-3)*{}} ;
  {\ar@{=>}_{\scriptstyle \alpha}(-7.5,3)*{};(-7.5,-3)*{}} ;
 \endxy
\]
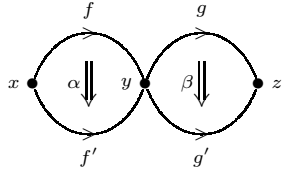
17.  Horizontal composition of 2-morphisms
\efig

\noindent One can show that vertical and horizontal composition satisfy an
`exchange identity'
\be    (\alpha \alpha') \tensor (\beta\beta') = (\alpha \tensor
\beta)(\alpha' \tensor \beta')  \label{exchange} \ee
making the diagram in Figure 18 define a unique 2-morphism.
This makes it quite convenient to define 2-morphisms
diagrammatically by `pasting' together diagrams, in a manner nicely
mimicking how one can paste together 2-manifolds with corners.
This is the basic sense in which 2-categories encode
2-dimensional topology.

\bfig
\[
 \xy 
  (-15,0)*{\bullet}="1";
  (0,0)*{\bullet}="2";
  (15,0)*{\bullet}="3";
 "1";"2" **\crv{(-12,9) & (-3,9)};
 "1";"2" **\crv{};
 "2";"3" **\crv{};
  "1";"2" **\crv{(-12,-9) & (-3,-9)};
   "2";"3" **\crv{(3,9) & (12,9)} ;
    "2";"3" **\crv{ (3,-9)& (12,-9) };
    (-7.5,6.75)*{\scriptstyle >};
    (7.5,6.75)*{\scriptstyle >};
    (-7.5,-6.75)*{\scriptstyle >};
    (-7.5,0)*{\scriptstyle >};
    (7.5,0)*{\scriptstyle >};
    (7.5,-6.75)*{\scriptstyle >};
  {\ar@{=>}_{\scriptstyle \beta'}(7.5,5)*{};(7.5,1.25)*{}} ;
  {\ar@{=>}_{\scriptstyle \beta}(7.5,-1.25)*{};(7.5,-5)*{}} ;
  {\ar@{=>}_{\scriptstyle \alpha'}(-7.5,5)*{};(-7.5,1.25)*{}} ;
  {\ar@{=>}_{\scriptstyle \alpha}(-7.5,-1.25)*{};(-7.5,-5)*{}} ;
 \endxy
\]
18.  Exchange identity
\efig

Just as the primordial example of a category is $\set$, the
primordial example of a 2-category is $\cat$.  We have already
discussed $\cat$ as a category in which the objects are
small categories and the morphisms are functors.
Actually, however, $\cat$ is a 2-category, in which given
functors $F,G \maps C \to D$, the 2-morphisms from $F$ to $G$
are the `natural transformations' $\alpha \maps F
\doublearrow G$.  Recall that such a thing
assigns to each object $x$ of $C$ a morphism $\alpha_x \maps
F(x)\to G(x)$,  in such a way that for any morphism $f \maps x
\to y$ in $C$, the consistency condition
\be \label{natural}
\begin{diagram}[F(x)]
\node{F(x)} \arrow{e,t}{\alpha_x} \arrow{s,l}{F(f)} \node{G(x)}
\arrow{s,r}{G(f)}       \\
\node{F(y)}  \arrow{e,t}{\alpha_y}    \node{G(y)}
\end{diagram}
\ee
holds.

While this example may seem abstract, it has a certain inherently
geometrical character.  A functor $F \maps C \to D$ can be
viewed as a diagram in $D$ shaped liked $C$, and a natural
transformation $\alpha \maps F \doublearrow G$ should then be
viewed as a prism in $D$ going between two such diagrams, as
shown in Figure 19.

\bfig
\[
 \xy 
    (-26,-2)*+{\bullet}="a";
    (-12,3)*+{\bullet}="b";
    (-16,-5)*+{\bullet}="c";
    {\ar@/^.25pc/ "a";"b"};
    {\ar@/_.25pc/"a";"c"};
    {\ar@/_.15pc/ "b";"c"};
    %START TOP
    (12,10)*+{\bullet}="ta";
    (26,15)*+{\bullet}="tb";
    (22,8)*+{\bullet}="tc";
    {\ar@/^.25pc/ "ta";"tb"};
    {\ar@/_.25pc/ "ta";"tc"};
    {\ar@/_.15pc/ "tb";"tc"};
        %START bottom
    (12,-10)*+{\bullet}="ba";
    (26,-5)*+{\bullet}="bb";
    (22,-12)*+{\bullet}="bc";
    {\ar@/^.25pc/@{.>} "ba";"bb"};
    {\ar@/_.25pc/ "ba";"bc"};
    {\ar@/_.15pc/ "bb";"bc"};
    %START CONNECT
    {\ar@/_.15pc/ "ta";"ba"};
    {\ar@/_.15pc/ "tb";"bb"};
    {\ar@/_.15pc/ "tc";"bc"};
    %START ARROWS AND LABELS
    (-12,6)*{}="DF";
    (6,10)*{}="CF";
    {\ar@/^.35pc/^F "DF";"CF"};
    (-12,-6)*{}="DF";
    (6,-10)*{}="CF";
    {\ar@/_.35pc/_G "DF";"CF"};
    (-2,5)*{}="DF";
    (-2,-5)*{}="CF";
    {\ar@{=>}^{\alpha} "DF";"CF"};
    (29,3)*{D};
    (-23,3)*{C};
 \endxy
\]
19.  Natural transformation $\alpha$ between
functors $F,G\maps C\to D$
\efig

\noindent The consistency condition in the definition of a
natural transformation says that the rectangular `vertical' faces
of this prism commute.

These remarks about 2-categories generalize considerably.   In an
$n$-category one has composition operations that allow one to
paste together $n$-morphisms according to a wide variety of
`pasting schemes', much as one can glue together $n$-manifolds
with corners \cite{Johnson}.  Moreover, the primordial example of
an $(n+1)$-category is $n\Cat$.  The reason is simply that
$n\Cat$ is closed, i.e., enriched over itself.    That is, in
addition to `$n$-functors' between $n$-categories and
`$n$-natural transformations' between these, there are higher
transformations between these which can be visualized using
higher-dimensional analogs of Figure 19.  Given two
$n$-categories $C$ and $D$, we thus obtain an $n$-category $\hom(C,D)$.

\section{Weakening}

One
profound difference between a set and a category is that elements
of a set are either {\it equal} or not, while objects
in a category can also be {\it isomorphic} in different ways (or
not at all).   Modern mathematics and physics takes advantage of
this insight in many ways.  For example, the fact that an object can
admit nontrivial automorphisms is precisely what yields the
notion of symmetry group.   However, it is primarily
category theorists who have followed through on this insight with
the philosophy of `weakening'.  As clearly enunciated
by Kapranov and Voevodsky \cite{KV}, this is based on
the principle that {\it ``In any
category it is unnatural and undesirable to speak about equality of
two objects.''}   Instead, it is better whenever possible to
speak in terms of isomorphisms between them.

For example, in the context of set
theory, algebraic structures are frequently defined using
equations.   These structures can often be generalized to the
context of category theory, but one has the choice of
generalizing them `strictly' --- keeping the equations as
equations --- or `weakly' --- replacing the equations by
specified isomorphisms.   When one opts to `weaken' a
definition in this way, one typically demands that the
isomorphisms themselves satisfy new equations, called `coherence
laws', in order to manipulate them with some of the same facility
as the original equations.  For example, a monoid is a set with a
product that is required among other things to satisfy the
equation $(xy)z = x(yz)$.   The categorical analog of a monoid is
a category with tensor product, or monoidal category.  Actually,
though, monoidal categories come in two versions: strict, where
the associativity of the tensor product is given by an equation:
\[
(x \tensor y) \tensor z = x \tensor (y \tensor z), \]
and `weak',
where instead there is a natural isomorphism, the `associator':
\[                A_{x,y,z}\maps (x \tensor y) \tensor z \to x
\tensor (y \tensor z).  \]
The associator allows one to rebracket
iterated tensor products, but to make sure that any two different
paths of rebracketings have  the same effect one must impose the
`pentagon identity', that
\[    \begin{diagram}[((x \tensor y)\tensor z)\tensor w]
\node{((x \tensor y) \tensor z)\tensor w}
\arrow{e,t}{A_{x\tensor y,z,w}}
\arrow{s,l}{A_{x,y,z}\tensor 1_w}
\node{(x\tensor y)\tensor(z\tensor w)}
\arrow{e,t}{A_{x,y,z\tensor w}}
\node{x \tensor (y\tensor(z \tensor w))}  \\
\node{(x \tensor (y \tensor z)) \tensor w}
\arrow[2]{e,t}{A_{x,y\tensor z,w}}
\node[2]{x \tensor ((y \tensor z)\tensor w)}
\arrow{n,r}{1_x \tensor A_{y,z,w}}
\end{diagram} \]
commutes.  Similarly, in a weak
monoidal category the equations $1x = x1 = x$ holding in a monoid
are replaced by isomorphisms satisfying coherence laws.

The monoidal categories that arise in nature, such as $n\Cob$ and
$\Vect$, are usually weak.   People frequently ignore this fact,
however (and the reader will note we did so in Section 1).
The justification for doing so is Mac Lane's theorem
\cite{Maclane} that any weak monoidal category is equivalent to a
strict one.   However, the sense of `equivalence' here is rather
subtle and itself intimately connected with weakening.
Following Kapranov and Voevodsky's principle, in
addition to weakening algebraic structures, one should also weaken
the sense in which maps between them preserve the structure.
For example, the strictest notion of a `monoidal functor'
between monoidal categories would require that it
preserve tensor products `on the nose'.  A weaker and often
more useful notion, however, requires merely that it preserve
tensor products {\it up to a natural isomorphism} compatible
with the associativity constraints.  It is this
weaker notion which plays a role in the definition of
`equivalence' of monoidal categories.

Mac Lane's theorem is an example of the `strictification' theorems in
higher- \break dimensional algebra.  These assert that any weakened
algebraic structure of a given sort is equivalent to some stricter
counterpart, in an appropriately weakened sense of `equivalence'.  They
simplify certain computations by allowing us to consider a special class
of cases without essential loss of generality.

However, in many situations weak notions are more general than
their strict counterparts in interesting ways.  Also,
there is a certain matter of choice involved in picking
coherence laws, and this can lead to different degrees of weakening.
For example, one
weakened categorical analog of a commutative monoid is a
symmetric monoidal category, the equation  $xy = yx$ having been
replaced by an isomorphism $B_{x,y} \maps x \tensor y \to y
\tensor x$ satisfying various coherence laws including $B_{y,x}B_{x,y} = 1_{x
\tensor y}$.  The stricter notion where commutativity remains
an equation is too narrow to be very interesting.  However, a still weaker
notion {\it is} very interesting, namely a `braided' monoidal
category, in which the coherence law $B_{y,x}B_{x,y} = 1_{x
\tensor y}$ is dropped.

As the name suggests, braided monoidal
categories are important in 3-dimensional topology
\cite{CP,FY,RT}.  A further `categorification' of the
notion of commutative monoid, namely a braided monoidal
2-category, appears to play a corresponding role in 4-dimensional
topology \cite{CS,KV}.   One goal of the $n$-categorical approach
to TQFTs is to systematically understand why weakened categorical
analogs of familiar algebraic structures are important in
topology.

Perhaps the most fundamental candidate for weakening is the
definition of $n$-category itself.   In the context of $n$-categories,
Kapranov and Voevodsky's principle indicates that it is
undesirable to speak of equality between two $k$-morphisms
when $k < n$; instead, one should speak in terms of $(k+1)$-isomorphisms
between them.  One can unfold the recursive
definition of strict $n$-category and obtain a
completely `explicit' definition in terms of operations
on $k$-morphisms, which are required to satisfy various
equations.  Each equation presents an opportunity for repeated
weakening.  For example, we can weaken an equation between
$k$-morphisms by replacing it with a natural $(k+1)$-isomorphism and
demanding that this new isomorphism satisfy coherence laws.
We can weaken further by replacing these coherence laws with
$(k+2)$-isomorphisms, which must satisfy their own coherence
laws, and so on.   This process becomes increasingly complex
with increasing $n$, and so far the definition of what might be called
a `weak $n$-category' has only been worked out for $n \le 3$.

In defining $n$-categories with $n = 0$ or $1$ --- i.e., in
defining sets or categories --- there is no opportunity for
weakening.  Weak 2-categories are usually known as `bicategories'
\cite{Benabou}, but there is a strictification theorem saying
that all of these are equivalent (or more precisely,
biequivalent) to strict 2-categories.  Weak 3-categories, or
`tricategories', have recently been developed by Gordon, Power,
and Street \cite{GPS}.  These are {\it not} all triequivalent to
strict 3-categories, but there {\it is} a strictification theorem
saying they are triequivalent to `semistrict 3-categories'.
These are categories enriched over $2\Cat$ thought of as a
monoidal category not with its Cartesian product, but  with a
weakened product similar to that defined by Gray \cite{Gray}.  In
this `semistrict' tensor product, eq.\ (\ref{cartprod}) is dropped,
and instead there is only a natural 2-isomorphism between the
left- and right-hand sides.   Topologically this is very natural,
since it means that the squares in Figure 20, rather than
commuting, are `filled in' with 2-isomorphisms.

\bfig
\[
 \xy 
    (-13,-2)*+{\bullet}="1";
    (-13,6)*+{\bullet}="2";
    (0,-2)*+{\bullet}="d";
    (3,1)*+{}="dx";
    (0,6)*+{\bullet}="a";
    (10,-2)*+{\bullet}="e";
    (13,1)*+{}="ex";
    (10,6)*+{\bullet}="b";
    (20,-2)*+{\bullet}="f";
    (20,6)*+{\bullet}="c";
    (0,15)*+{\bullet}="4";
    (10,15)*+{\bullet}="5";
    (20,15)*+{\bullet}="6";
    (10,18)*{C};
    (-18,2.5)*{D};
    (29,2.5)*{C \tensor D};
    {\ar@/^.25pc/ "2";"1"};
    {\ar@/^.25pc/ "a";"b"};
    {\ar@/^.25pc/ "b";"c"};
    {\ar@/^.25pc/ "d";"e"};
    {\ar@/^.25pc/ "e";"f"};
    {\ar@/^.25pc/ "a";"d"};
    {\ar@/^.25pc/ "b";"e"};
    {\ar@/^.25pc/ "c";"f"};
    {\ar@/^.25pc/ "4";"5"};
    {\ar@/^.25pc/ "5";"6"};
    {\ar@{<=>} "b";"dx"};
    {\ar@{<=>} "c";"ex"};
 \endxy
\]
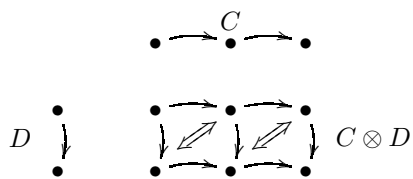
20.  The semistrict tensor product of 2-categories
\efig

To advance further in $n$-category theory, it is urgent to define
`weak $n$-categories' for all $n$.  It is clear that new ideas are
needed to do so without a combinatorial explosion, since already the
explicit definition of a tricategory takes 6 pages, and that of a
triequivalence 13 pages!  However, the potential payoffs of a good
theory of weak $n$-categories should encourage us to persevere.

Having completed our brief survey of $n$-category theory, let us
return to topological quantum field theory.   In what follows, we
propose answers to the basic questions: {\it Of which
$n$-category are $n$-dimensional extended TQFTs representations?}
and {\it In what sense is an $n$-dimensional extended TQFT a
representation of this $n$-category?}  Our answers are inevitably
somewhat vague, except for low $n$, since they rely on notions
from the theory of weak $n$-categories.   Nonetheless, we hope
they will serve as a guide for future research.

\section{Suspension}

To begin, it is useful to consider an issue that might at first
seem of purely formal interest.  Suppose we have an
$(n+1)$-category $C$ with only one object $x$.    We can
regard $C$ as an $n$-category $\tilde C$ by re-indexing: the
objects of $\tilde C$ are the morphisms of $C$, the morphisms of
$\tilde C$ are the 2-morphisms of $C$, and so on.  However, the
$n$-categories we obtain this way will have extra structure.  For
example, since the objects of $\tilde C$ were really morphisms in
$C$ from $x$ to itself, we can multiply (i.e., compose) them.
We have already seen the simplest example of this phenomenon in
Section 1: if $C$ is a category with a single object, $\tilde C$
is a monoid.   If instead we start with a strict (resp.\ weak)
2-category with a single object, we obtain a strict (resp.\ weak)
monoidal category!   Similarly, starting with strict, semistrict,
or weak 3-categories with only one object, we obtain
corresponding sorts of monoidal 2-categories, i.e., 2-categories
having tensor products of objects, morphisms, and 2-morphisms
\cite{GPS,KV}.

We can iterate this process, and construct from an
$(n+k)$-category $C$ with only one object, one  morphism, and so
on up to one $(k-1)$-morphism, an $n$-category $\tilde C$ whose
$j$-morphisms are the $(j+k)$-morphisms of $C$.  In doing so we
obtain a particular sort of $n$-category with extra structure and
properties, which we call a `$k$-tuply monoidal' $n$-category.  In
Figure 21 we tabulate our best guesses concerning  $k$-tuply
monoidal $n$-categories.  Ultimately we expect a table along
these lines for {\it weak} $k$-tuply monoidal $n$-categories.
For the moment, however, we work with `semistrict' ones,
which have already been defined in a few cases where the weak ones have not.
The idea is that strictification theorems are either known or
expected saying that all weak $k$-tuply monoidal
$n$-categories are equivalent (in a suitable sense) to these
semistrict ones.

More precisely, for the $n = 0$ and $n = 1$ columns we define the
semistrict notions to be the same as the strict ones.  In the $n
= 2$ column, we define semistrict 2-categories to be strict ones,
while semistrict monoidal and braided monoidal 2-categories have
been defined by Kapranov and Voevodsky \cite{KV}.  Semistrict
weakly and strongly involutory monoidal 2-categories have been
discussed by Breen \cite{Breen}.  Semistrict 3-categories,
mentioned in the previous section, have been studied by Gordon,
Power and Street \cite{GPS} and Leroy \cite{Leroy}.

\vskip 0.5em
\begin{center}
{\small
\begin{tabular}{|c|c|c|c|}  \hline
         & $n = 0$   & $n = 1$    & $n = 2$          \\     \hline
$k = 0$  & sets      & categories & 2-categories     \\     \hline
$k = 1$  & monoids   & monoidal   & monoidal         \\
         &           & categories & 2-categories     \\     \hline
$k = 2$  &commutative& braided    & braided          \\
         & monoids   & monoidal   & monoidal         \\
         &           & categories & 2-categories     \\     \hline
$k = 3$  &`'         & symmetric  & weakly involutory \\
         &           & monoidal   & monoidal         \\
         &           & categories & 2-categories     \\     \hline
$k = 4$  &`'         & `'         &strongly involutory\\
         &           &            & monoidal         \\
         &           &            & 2-categories     \\     \hline
$k = 5$  &`'         &`'          & `'               \\
         &           &            &                  \\
         &           &            &                  \\      \hline
\end{tabular}} \vskip 1em
21.  Semistrict $k$-tuply monoidal $n$-categories
\end{center}
\vskip 0.5em

There are many interesting patterns to be seen in this table.
First, it is clear that many of the concepts already discussed
appear in this table, together with some new ones.
Second, as we proceed down any column of this table, the
$n$-categories in question first gain additional
structures, which then acquire additional properties of an
`abelian' nature.   This process appears in its most rudimentary
form in the first column.   We have already seen that a category
with only one object $x$ is essentially the same as the monoid
$\hom(x,x)$.  Why does a 2-category $C$ with only one object $x$ and
one morphism $1_x$ give a {\it commutative} monoid $\hom(1_x,1_x)$?

The argument goes back at least to Eckmann and Hilton \cite{EH}.
The elements of $\hom(1_x,1_x)$ are the 2-morphisms of $C$, and
as described in our brief review of 2-categories, we can compose
$\alpha,\beta\maps 1_x \doublearrow 1_x$ either vertically or
horizontally to obtain a new 2-morphism from $1_x$ to $1_x$.
We write the vertical composite as $\alpha\beta$ and the
horizontal composite as $\alpha\tensor \beta$.
Writing simply $1$ for $1_{1_x}$, with a little work one
can check that $1 \tensor \alpha = \alpha \tensor 1 = \alpha$, so
that $1$ is the identity for both vertical and horizontal
composition.  One also has the exchange identity
$(\alpha \alpha') \tensor (\beta\beta') = (\alpha \tensor
\beta)(\alpha' \tensor \beta')$ from eq.\ (\ref{exchange}).
These two facts let us perform the remarkable computation:
\ba   \alpha\tensor\beta  &=& (1 \alpha) \tensor (\beta 1)\label{eh}\\
	&=& (1 \tensor \beta)(\alpha \tensor 1) \nonumber\\
     &=& \beta\alpha \nonumber\\
	&=& (\beta \tensor 1)(1 \tensor \alpha)  \nonumber\\
     &=& (\beta 1) \tensor (1 \alpha) \nonumber\\
     &=& \beta\tensor\alpha, \nonumber\ea
so vertical and horizontal composition are equal and
$\hom(1_x,1_x)$ is a commutative monoid.  Conversely, one can
show that any commutative monoid can be thought of as the
2-morphisms in a 2-category with one object and one morphism.

When we consider a semistrict 3-category with one object $x$, one
morphism $1_x$, and one 2-morphism $1_{1_x}$, it turns out that
these are again essentially just commutative monoids.   The same
appears to be true for semistrict 4-categories with only one
3-morphism.  While semistrict 4-categories are not understood in
general, it seems that those with only one morphism are,
these being braided monoidal 2-categories, and one can check that
of these, those with only one 3-morphism are commutative monoids.
Of course this argument is somewhat circular, since it assumes we
understand one column of the table in order to check another, but
it serves as a interesting cross-check.   In any event, it
appears that the $n = 0$ column {\it stabilizes} after two steps.

The same sort of process is at work in the next two columns, in
increasingly sophisticated incarnations.  The basic pattern is
shown in Figure 22.

\vskip 0.5em
\begin{center}
{\small
\begin{tabular}{|c|c|c|c|}  \hline
         & $n = 0$   & $n = 1$      & $n = 2$          \\     \hline
$k = 0$  &    sets   & categories   & 2-categories     \\     \hline
$k = 1$  &    $xy$   &$x \tensor y$ & $x \tensor y$    \\     \hline
$k = 2$  & $xy = yx$ & $B_{x,y}\maps x\tensor y \to y\tensor x$
 & $B_{x,y}\maps x\tensor y \to y\tensor x$            \\      \hline
$k = 3$  &`'         & $B_{x,y} = B_{y,x}^{-1}$  & $I_{x,y} \maps
B_{x,y} \doublearrow B_{y,x}^{-1}$                       \\      \hline
$k = 4$  &`'         & `'         & $I_{x,y} =
(I_{y,x}^{-1})^{-1}_{hor}$                              \\     \hline
\end{tabular}} \vskip 1em
22.  Semistrict $k$-tuply monoidal $n$-categories:
structure and properties
\end{center}
\vskip 0.5em

\noindent In the $n = 0$ column we began with sets, which then
acquired a product, which then satisfied the commutativity
equation $xy = yx$.  In the $n = 1$ we begin with categories,
which permit a more nuanced version of the process: first they
acquire a product, then they acquire an isomorphism $B_{x,y}
\maps x\tensor y \to y \tensor x$, taking the place of the
commutativity equation.  Finally, the braiding is required to
satisfy an equation of its own, the symmetry equation
$B_{x,y} = B_{y,x}^{-1}$.  Note that what was a property
(commutativity) has become structure (the braiding), which then
acquires an analogous property of its own (symmetry).

In the $n = 2$ column a still more subtle version of the
`abelianization' process occurs.  We begin with 2-categories.
These first acquire a product,  then a braiding isomorphism, and
then a 2-isomorphism $I_{x,y} \maps B_{x,y} \doublearrow
B_{y,x}^{-1}$, the `involutor', taking the place of the symmetry
equation.  In the last step, the involutor satisfies an equation
of its own,  $I_{x,y} = (I_{y,x}^{-1})^{-1}_{hor}$, meaning
that the horizontal composite of $I_{x,y}$ and
$I_{y,x}^{-1}$ is the identity.   One can also think of both
sides of the equation as 2-isomorphisms from $B_{x,y}$ to
$B_{y,x}^{-1}$.  We give a topological interpretation of this
equation in Section 7.

The pattern here is evident, at least in outline, and it is
tempting to predict that it continues for higher $n$.  In
particular, we can guess that each column will take one step
longer to stabilize.   A bit more precisely,
let $n\cat_k$ denote the category of $k$-tuply monoidal
weak $n$-categories.   There should be a forgetful functor
\[        F\maps n\cat_k \to n\cat_{k-1} , \]
and a corresponding `reverse' functor, technically a left adjoint
\[       S\maps n\cat_{k-1} \to n\cat_k  ,\]
which we shall call `suspension' (for reasons to become clear
shortly).   For example, when we repeatedly suspend a set $C$, we
obtain first the free monoid on $C$, and  then the abelianization
thereof.  Similarly, when we repeatedly suspend a category $C$, we
obtain first the free monoidal category on $C$, then the free
braided monoidal category on $C$, and then the symmetrization
thereof.   We propose the:
\vskip 1em\noindent
{\bf Stabilization Hypothesis.}  {\it After suspending a
weak $n$-category $n+2$ times, further suspensions have no
essential effect.  More precisely, the suspension functor $S \maps
n\cat_k \to n\cat_{k+1}$ is an equivalence of
categories for $k \ge n+2$.}
\vskip 1em
\noindent We could also embellish this hypothesis by considering
$n\cat_{k}$ not merely as a category but as $(n+k+1)$-category,
and using not equivalence but some notion of
`$(n+k+1)$-equivalence'.

One can regard the above hypothesis, and those to follow,
either as a conjecture pending a general definition of
`weak $n$-category', or as a feature one might desire of
such a definition.  Apart from the already given algebraic
evidence for the stabilization hypothesis, there is quite a bit
of indirect topological evidence from homotopy theory and
the theory of tangles.  The latter leads us back to our goal:
understanding topological quantum field theory in $n$-categorical
terms.

\section{Homotopy Theory}

Modern higher-dimensional algebra has it roots in the dream of
finding a natural and convenient {\it completely algebraic}
description of the homotopy type of a topological space.  The
prototype here is the fundamental groupoid; given a space $X$,
this is the category $\Pi_1(X)$ whose objects are the points of
$X$ and whose morphisms are the homotopy classes of paths (with
fixed endpoints).  This is a groupoid, meaning that every
morphism has an inverse, given by reversing a path from $x$ to
$y$ to obtain a path from $y$ to $x$.  In fact, $\Pi_1$ gives an
equivalence between the category of groupoids and the category of
`homotopy 1-types', where, roughly speaking, two spaces define the
same homotopy $n$-type if they are equivalent as far as concerns
homotopy classes of maps from $n$-dimensional CW complexes into
them.

The goal of generalizing the fundamental groupoid to higher dimensions
has led to a variety of schemes.  One of the most popular involves Kan
complexes \cite{May}, which model a space by an algebraic analog of a
simplicial complex.  Alternative approaches based on cubes have been
developed by Brown, Higgins, Loday and others \cite{Brown}.  Indeed, it
was in this context that Brown first coined the term `higher-dimensional
algebra'.

Here, however, we restrict our attention to $n$-categorical
approaches.  Ever since Grothendieck's famous 600-page letter to Quillen
\cite{Gro}, it has been tempting to associate to a space $X$ a
`fundamental $n$-groupoid' $\Pi_n(X)$, some sort of $n$-category
whose objects are points, whose morphisms are paths, whose
2-morphisms are paths between paths, and so on up to the
$n$-morphisms, which are homotopy classes of $n$-fold paths.
This amounts to taking the imagery of $n$-category theory quite
literally, as in Figure 23, where we show a typical 2-morphism.

\bfig
\[
 \xy 
  (-15,0)*{\bullet}="1"; 
  (0,0)*{\bullet}="2";
 "1";"2" **\crv{(-12,9) & (-3,9)};
  "1";"2" **\crv{(-12,-9) & (-3,-9)};
    (-7.5,6.75)*{\scriptstyle >}+(0,3)*{ };
    (-7.5,-6.75)*{\scriptstyle >}+(0,-3)*{ };
    (-17.5,0)*{ };
  (2.5,0)*{ };
  {\ar@{=>}_{ }(-7.5,3)*{};(-7.5,-3)*{}} ;
 \endxy
\]
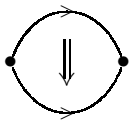
23.  A 2-morphism in the fundamental $n$-groupoid
\efig

Roughly speaking, an $n$-groupoid should be some sort of
$n$-category in which all $k$-morphisms ($k \ge 1$) have
inverses, at least weakly.   There are good reasons to want to
use weak $n$-categories here.  In the fundamental groupoid,
composition is associative, since the morphisms are merely
homotopy classes of paths.  In the fundamental 2-groupoid,
however, composition of paths $f \maps [0,1] \to X$ is not
strictly associative, but only up to a homotopy, the associator,
which performs the reparametrization of $[0,1]$ shown in Figure
24.

\bfig
\[
 \xy 
    (-15,7)*{}="TL";
    (-15,-7)*{}="BL";
    (15,7)*{}="TR";
    (15,-7)*{}="BR";
    "TL";"BL" **\dir{-};
    "TL";"TR" **\dir{-};
    "TR";"BR" **\dir{-};
    "BL";"BR" **\dir{-};
    (-10,7)*{};(0,-7)*{} **\dir{-};
    (0,7)*{};(10,-7)*{} **\dir{-};
    (-10,10)*{(fg)};
    (7,10)*{h};
    (10,-10)*{(gh)};
    (-7,-10)*{f};
 \endxy
 \]
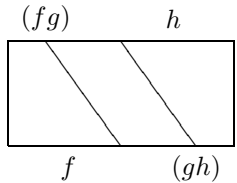
24.  The associator in homotopy theory
\efig

\noindent In the fundamental 2-groupoid, the associator satisfies
the pentagon identity `on the nose', but in higher fundamental
$n$-groupoids it does so only up to a homotopy, which in turn
satisfies a coherence condition up to homotopy, and so on.  In
fact, the whole tower of these `higher associativity laws' was
worked out by Stasheff \cite{Stasheff} in 1963, and have an
appealing geometrical description as faces of the
`associahedron'.  For example, the pentagon, being a 2-morphism,
is a 2-dimensional face.  One expects these higher associativity
laws to play a crucial part in the definition of weak
$n$-categories, as indeed they do in the cases understood so far
($n \le 3$).   Similar remarks hold for the identity law
\cite{GPS} and for the inverses.

Of course, one expects strictification theorems saying that `weak
$n$-groupoids' are all `$n$-equivalent' to some better-behaved class of
$n$-groupoids, implying that the latter are sufficient for homotopy
theory.  There are, in fact, two distinct strands of progress along
these lines.  First, the category of homotopy $2$-types has been shown
equivalent to a category whose objects are strict 2-categories having
strict inverses for all $k$-morphisms \cite{BS,MW}.  Moreover, the
category of homotopy 3-types has been shown equivalent to a category
whose objects are semistrict 3-categories having strict inverses
\cite{JT,Leroy}.  This naturally suggests the possibility that homotopy
$n$-types might be equivalent to some sort of semistrict $n$-categories
having strict inverses.  Second, for all $n$ the category of homotopy
$n$-types has been shown equivalent to a category whose objects are
strict $n$-categories having a particular sort of `weak inverses'
\cite{KVinfinity}.  We see here a kind of tradeoff that would be nice to
understand better.

In any event, while the correspondence between homotopy $n$-types
and `weak $n$-groupoids' is still incompletely understood,
it is already very valuable, since it sets up an analogy between
topological spaces and $n$-categories that lets us import
techniques and insights from topology into higher-dimensional
algebra. In what follows we use this analogy to shed some light
on the stabilization hypothesis of the previous section.

First, the topologically minded reader might already have noticed
that the $n = 0$ column of Figure 21 --- sets, monoids, and
commutative monoids thereafter --- is familiar from homotopy
theory.  Typically in homotopy theory one works with spaces with
basepoint, and defines $\pi_k(X)$ to be the set of homotopy
classes of based maps from $S^n$ to $X$.  For $k = 0$, $\pi_k(X)$
is indeed just a set, while for $k = 1$ it is a group, and for $k
\ge 2$ it is an abelian group.  These facts can be seen as {\it
consequences} of Figure 21 together with the correspondence
between homotopy $n$-types and weak $n$-groupoids, but
historically, of course, they were discovered first.

In particular, the Eckmann-Hilton argument given in eq.\
(\ref{eh}) is just the algebraic distillation of a very topological
proof that $\pi_2$ is abelian.  We illustrate the key steps of
eq.\ (\ref{eh}) in Figure 25.  Here each rectangle is
labelled with a map from the rectangle to some space $X$,
mapping the boundary to the basepoint of $X$.  Such a map can be
thought of as a based map from $S^2$ to $X$, defining an element
of $\pi_2(X)$.   In Figure 25, $\alpha$ and $\beta$ represent arbitrary
based maps from $S^2$ to $X$, while $1$ represents the trivial
map sending all of $S^2$ to the basepoint of $X$.  Such maps can
be either vertically or horizontally composed by juxtaposition.
and the figure shows successive frames in a movie of a homotopy
from $\alpha \tensor \beta$ to $\beta \tensor \alpha$, proving
that $\pi_2(X)$ is abelian.   The same sort of argument shows that
$\pi_k(X)$ is abelian for all $k \ge 2$.

\bfig
\[ 
 \xy
 (-6,6)*{}="TL";
 (-6,-6)*{}="BL";
 (6,6)*{}="TR";
 (6,-6)*{}="BR";
 (0,6)*{}="MT";
 (0,-6)*{}="MB";
 "TL";"TR" **\dir{-};
 "TL";"BL" **\dir{-};
 "TR";"BR" **\dir{-};
 "BL";"BR" **\dir{-};
 "MT";"MB" **\dir{-};
 (3,0)*{\beta};
 (-3,0)*{\alpha};
 (0,-9)*{\scriptstyle \alpha \tensor \beta};
 \endxy
\qquad
 \xy
 (-6,6)*{}="TL";
 (-6,-6)*{}="BL";
 (6,6)*{}="TR";
 (6,-6)*{}="BR";
 (0,6)*{}="MT";
 (0,-6)*{}="MB";
 (-6,0)*{}="LM";
 (6,0)*{}="LR";
 "TL";"TR" **\dir{-};
 "TL";"BL" **\dir{-};
 "TR";"BR" **\dir{-};
 "BL";"BR" **\dir{-};
 "MT";"MB" **\dir{-};
  "LM";"LR" **\dir{-};
 (3,-3)*{\beta};
 (-3,3)*{\alpha};
 (-3,-3)*{1};
 (3,3)*{1};
 (0,-9)*{\scriptstyle (1 \tensor \beta)(\alpha \tensor 1)};
 \endxy
 \qquad
  \xy
 (-6,6)*{}="TL";
 (-6,-6)*{}="BL";
 (6,6)*{}="TR";
 (6,-6)*{}="BR";
 (-6,0)*{}="LM";
 (6,0)*{}="LR";
 "TL";"TR" **\dir{-};
 "TL";"BL" **\dir{-};
 "TR";"BR" **\dir{-};
 "BL";"BR" **\dir{-};
  "LM";"LR" **\dir{-};
 (0,-3)*{\beta};
 (0,3)*{\alpha};
 (0,-9)*{\scriptstyle \beta\alpha};
 \endxy
 \qquad
  \xy
 (-6,6)*{}="TL";
 (-6,-6)*{}="BL";
 (6,6)*{}="TR";
 (6,-6)*{}="BR";
 (0,6)*{}="MT";
 (0,-6)*{}="MB";
 (-6,0)*{}="LM";
 (6,0)*{}="LR";
 "TL";"TR" **\dir{-};
 "TL";"BL" **\dir{-};
 "TR";"BR" **\dir{-};
 "BL";"BR" **\dir{-};
 "MT";"MB" **\dir{-};
  "LM";"LR" **\dir{-};
 (-3,-3)*{\beta};
 (3,3)*{\alpha};
 (3,-3)*{1};
 (-3,3)*{1};
 (0,-9)*{\scriptstyle (\beta  \tensor 1)(1 \tensor \alpha)};
 \endxy
 \qquad
  \xy
 (-6,6)*{}="TL";
 (-6,-6)*{}="BL";
 (6,6)*{}="TR";
 (6,-6)*{}="BR";
 (0,6)*{}="MT";
 (0,-6)*{}="MB";
 "TL";"TR" **\dir{-};
 "TL";"BL" **\dir{-};
 "TR";"BR" **\dir{-};
 "BL";"BR" **\dir{-};
 "MT";"MB" **\dir{-};
 (-3,0)*{\beta};
 (3,0)*{\alpha};
 (0,-9)*{\scriptstyle  \beta \tensor \alpha };
 \endxy
\]
25.  The Eckmann-Hilton argument
\efig

Figure 25 also clarifies how, when we proceed from the $n = 0$ to
the $n = 1$ column in Figure 21, the commutativity equation is
weakened to a braiding isomorphism.  In $\pi_2(X)$, homotopic
maps from  $S^2$ to $X$ are decreed to be {\it equal}.  In a more
refined context, however, we could regard the homotopy as an {\it
isomorphism}.  In Figure 26 we depict the homotopy from
$\alpha\tensor\beta$ to $\beta\tensor\alpha$ given by the
Eckmann-Hilton argument as a map from the cube to $X$, whose
horizontal slices give the frames of the movie shown in Figure
25. (Here we have compressed $\alpha$ and $\beta$ to small
discs for clarity; everything outside these discs is mapped
to the basepoint of $X$.)  One can see that this homotopy is
precisely a braiding!

\bfig
\[ 
 \xy 0;/r.16pc/:
 %*********************BUILDS BOX*****************************************************
 (-14,10)*{}="TL";
 (14,10)*{}="TR";
 (14,-10)*{}="BR";
 (-14,-10)*{}="BL";
 (-6,20)*{}="xTL";
 (22,20)*{}="xTR";
 (22,0)*{}="xBR";
 (-6,0)*{}="xBL";
 (0,10)*{}="M1";
 (8,20)*{}="M2";
 (0,-10)*{}="MB1";
 (8,0)*{}="MB2";
 %*********************adjustments to box*************
 (5,0)*{}="LeftBackString";
 (2,0)*{}="RightFrontString";
 (-2.5,0)*{}="LeftFrontString";
 (6.75,-1.25)*{}="BottomBackString";
 (10,0)*{}="RightBackString";
%*********************************************************************
     "TL";"TR" **\dir{-};
     "TR";"BR" **\dir{-};
     "BR";"BL" **\dir{-};
     "BL";"TL" **\dir{-};
     "xTL";"xTR" **\dir{-};
     "xTR";"xBR" **\dir{-};
     "xBR";"RightBackString" **\dir{.};
     "xBL";"LeftFrontString" **\dir{.};
     "LeftBackString";"RightFrontString" **\dir{.};
     "TL";"xTL" **\dir{-};
     "TR";"xTR" **\dir{-};
     "BL";"xBL" **\dir{.};
     "BR";"xBR" **\dir{-};
     "M1";"M2" **\dir{-};
     "MB1";"BottomBackString" **\dir{.};
     "xTL";"xBL" **\dir{.};
 %*********************END BUILD BOX***************************************************
 (6,7.5)*\ellipse(2,1){-};
 (-1,-2.5)*\ellipse(2,1){-};
 (-2,7.5)*\ellipse(2,1){-};
 (5,-2.5)*\ellipse(2,1){-};
 (10,15)*{}="s1";
 (-4,-5)*{}="s2";
 "s1";"s2" **\crv{(10,5) & (-4,5)};
    \POS?(.67)*{\hole}="z2"; \POS?(.5)*{\hole}="z1";
 (14,15)*{}="s1";
 (0,-5)*{}="s2";
 "s1";"s2" **\crv{(14,5) & (0,5)}; \POS?(.71)*{}="y2"; \POS?(.61)*{}="y1";
 (-6,15)*{}="LL";
 (-2,15)*{}="LR";
 "LL";"z2" **\crv{(-6,4) & (0,5)};
 "LR";"z1" **\crv{(-2,4) & (4,5)};
 (8,-5)*{}="BL";
 (12,-5)*{}="BR";
 "y2";"BL" **\crv{(7.5,0)};
 "y1";"BR" **\crv{(10,4) & (12,-5)};
 (13,18)*{\scriptscriptstyle \beta};
  (-5,-8)*{\scriptscriptstyle \beta};
   (10,-8)*{\scriptscriptstyle \alpha};
    (-3,18)*{\scriptscriptstyle \alpha};
 \endxy
\]
26.  Braiding $B_{\alpha,\beta} \maps \alpha\tensor\beta \to
\beta\tensor\alpha$
\efig

In
particular, an alternate version of the Eckmann-Hilton argument:
\ban   \alpha\tensor\beta  &=& (\alpha 1)\tensor (1 \beta)\\
     &=& (\alpha \tensor 1)(1 \tensor \beta) \\
     &=& \alpha\beta \\
     &=& (1 \tensor \alpha)(\beta \tensor 1)\\
     &=& (1\beta) \tensor (\alpha 1) \nonumber\\
     &=& \beta\tensor\alpha, \ean
gives another homotopy from $\alpha\tensor\beta$ to
$\beta\tensor\alpha$, corresponding to $B_{\beta,\alpha}^{-1}$, as shown
in Figure 27.

\bfig
\[ 
 \xy 0;/r.16pc/:
 %*********************BUILDS BOX*****************************************************
 (-22,10)*{}="TL";
 (6,10)*{}="TR";
 (6,-10)*{}="BR";
 (-22,-10)*{}="BL";
 (-14,20)*{}="xTL";
 (14,20)*{}="xTR";
 (14,0)*{}="xBR";
 (-14,0)*{}="xBL";
 (-8,10)*{}="M1";
 (0,20)*{}="M2";
 (-8,-10)*{}="MB1";
 (0,0)*{}="MB2";
 %*********************adjustments to box*************
 (-3,0)*{}="LeftBackString";
 (-6,0)*{}="RightFrontString";
 (-10.5,0)*{}="LeftFrontString";
 (-1.25,-1.25)*{}="BottomBackString";
 (2,0)*{}="RightBackString";
%*********************************************************************
     "TL";"TR" **\dir{-};
     "TR";"BR" **\dir{-};
     "BR";"BL" **\dir{-};
     "BL";"TL" **\dir{-};
     "xTL";"xBL" **\dir{.};
     "xTL";"xTR" **\dir{-};
     "xTR";"xBR" **\dir{-};
     "xBR";"RightBackString" **\dir{.};
     "xBL";"LeftFrontString" **\dir{.};
     "LeftBackString";"RightFrontString" **\dir{.};
     "TL";"xTL" **\dir{-};
     "TR";"xTR" **\dir{-};
     "BL";"xBL" **\dir{.};
     "BR";"xBR" **\dir{-};
     "M1";"M2" **\dir{-};
     "MB1";"BottomBackString" **\dir{.};
 %*********************END BUILD BOX***************************************************
 (-6,7.5)*\ellipse(2,1){-};
 (1,-2.5)*\ellipse(2,1){-};
 (2,7.5)*\ellipse(2,1){-};
 (-5,-2.5)*\ellipse(2,1){-};
 (-10,15)*{}="s1";
 (4,-5)*{}="s2";
 "s1";"s2" **\crv{(-10,5) & (4,5)};
    \POS?(.67)*{\hole}="z2"; \POS?(.5)*{\hole}="z1";
 (-14,15)*{}="s1";
 (0,-5)*{}="s2";
 "s1";"s2" **\crv{(-14,5) & (0,5)}; \POS?(.71)*{}="y2"; \POS?(.61)*{}="y1";
 (6,15)*{}="LL";
 (2,15)*{}="LR";
 "LL";"z2" **\crv{(6,4) & (0,5)};
 "LR";"z1" **\crv{(2,4) & (-4,5)};
 (-8,-5)*{}="BL";
 (-12,-5)*{}="BR";
 "y2";"BL" **\crv{(-7.5,0)};
 "y1";"BR" **\crv{(-10,4) & (-12,-5)};
 (-12,18)*{\scriptscriptstyle \alpha};
  (2,-8)*{\scriptscriptstyle \alpha};
   (-12,-8)*{\scriptscriptstyle \beta};
    (5,18)*{\scriptscriptstyle \beta};
 \endxy
\]
27.  $B_{\beta,\alpha}^{-1} \maps \alpha\tensor\beta \to
\beta\tensor\alpha$
\efig

\noindent  These two homotopies are not generally homotopic to
each other, corresponding to the fact that generally $B_{\alpha,\beta} \ne
B_{\beta,\alpha}^{-1}$ in a braided monoidal category.  This is a
good example of how different {\it proofs} of the same equation
may, upon weakening, give rise to distinct isomorphisms.

Algebraically, the point here is that the Eckmann-Hilton
argument relies on the interplay between horizontal and vertical
composition of 2-morphisms in 2-categories. For the definition of
horizontal composition given in Figure 17 to be unambiguous, we
needed the commuting square condition, eq.\ (\ref{cartprod}),
in the definition of the Cartesian product of categories.  In
the context of semistrict 3-categories this equation is weakened
to an isomorphism, as shown in Figure 20, so Figure 17 gives two
distinct notions of horizontal composition of 2-morphisms,
related by a 3-isomorphism.  Arbitrarily picking one of these,
and recapitulating the Eckmann-Hilton argument in this context,
one finds that instead of an equation between $\alpha \tensor
\beta$ and $\beta \tensor \alpha$, one has a 3-isomorphism, the
braiding.

Lest the reader think we have drifted hopelessly far from physics
by now, we should note that elements of $\pi_2(X)$ correspond to
`topological solitons' in a nonlinear sigma model with target
space $X$, for a spacetime of dimension 3.  In this context,
Figures 26 and 27 show the worldlines of such topological
solitons, and the fact that the two pictures cannot be deformed
into each other is why the statistics of such solitons is
described using representations of the  braid group \cite{BOR}.
In a spacetime of dimension 4 or more, the analogous pictures
{\it can} be deformed into each other, since there is enough room
to pass the two strands across each other.   Topologically, this
means that the two homotopies from  $\alpha \tensor \beta$ to
$\beta \tensor \alpha$ are themselves homotopic when
$\alpha,\beta$ represent based maps from $S^k$ to $X$ for $k \ge
3$.  Algebraically, this corresponds to moving down in Figure 21
from braided monoidal categories to symmetric monoidal
categories, where one has $B_{\alpha,\beta} =
B_{\beta,\alpha}^{-1}$.  Physically, this is why the statistics of
topological solitons in spacetimes of dimension 4 or more is
described using the symmetric group.

Having discussed the $n = 0$ and $n = 1$ columns of Figure 21
from the viewpoint of homotopy theory, we could proceed to the $n
= 2$ column, but instead let us explain the general pattern and
what it has to do with the stabilization hypothesis.  In fact,
the `suspension' operation on $n$-categories is closely modelled
after homotopy theory.  In topology, one obtains the `suspension'
$SX$ of a space $X$ with basepoint $\ast$ as a quotient space of
$X \times [0,1]$ in which one collapses all the points of the
form $(x,0)$, $(x,1)$, and $(\ast,t)$ to a single point.  We can
draw this as in Figure 28, with the proviso that all the points
on the dotted line are identified.

\bfig
\[
 \xy 
  (0,0)*\ellipse(6,2.5){-};
  (-2,-2.35)*{\ast}="S";
  (-9,0)*{X};
 \endxy
\qquad \qquad \qquad
 \xy
  (0,0)*\ellipse(8,2.5){-};
  (-2,-2.35)*{\ast}="S";
  (-7.75,0.6)*{}="LT";
  (-7.75,-0.6)*{}="LB";
  (7.75,0.6)*{}="RT";
  (7.75,-0.6)*{}="RB";
  (0,10)*{\bullet}="T";
  (0,-10)*{\bullet}="B";
  "T";"LT" **\dir{-};
  "T";"RT" **\dir{-};
  "B";"LB" **\dir{-};
  "RB";"B" **\dir{-};
  "T";"S" **\dir{.};
  "S";"B" **\dir{.};
  (-12,0)*{SX};
 \endxy
\]
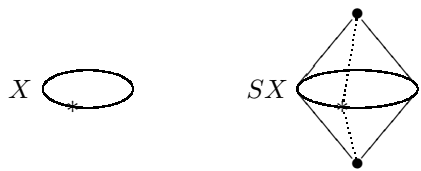
28.  Space $X$ and its suspension $SX$
\efig

\noindent In fact, suspension is a functor, so a map
$f \maps X \to Y$ gives rise to a map $Sf \maps SX \to SY$, and
one obtains thereby a sequence
\[ \begin{diagram}[(S^2X,S^2Y)]  \node{[X,Y]} \arrow{e,t}{S}
\arrow{e,t}{S} \node{[SX,SY]} \arrow{e,t}{S} \node{[S^2X,S^2Y]}
\arrow{e,t}{S} \node{\cdots }
\end{diagram} \]
where $[X,Y]$ denotes the set of homotopy classes of maps from
$X$ to $Y$.  Moreover, if $X$ is a CW complex of dimension $n$,
this sequence stabilizes after $n+2$ steps, i.e., the map
\[    S \maps [S^k X,S^k Y] \to [S^{k+1}X, S^{k+1}Y] \]
is an isomorphism for $k \ge n + 2$.  This theorem is the basis
of stable homotopy theory, a subject with close ties to
higher-dimensional algebra \cite{Adams}.   Given the conjectured
relation between homotopy $n$-types and weak $n$-groupoids,
one expects this theorem to translate into a proof of our
stabilization hypothesis in the special case of weak
$n$-groupoids.   Roughly speaking, the idea is that when we
suspend $X$, an $k$-morphism in the fundamental $n$-groupoid of $X$
gives a $(k+1)$-morphism in the fundamental $(n+1)$-groupoid of
$SX$ in a manner analogous to the algebraic notion of suspension
described in Section 5.  In Figure 29, for example, we show how a point
of $X$ gives a loop in $SX$, and how a path in $X$ gives a path
between loops in $SX$.

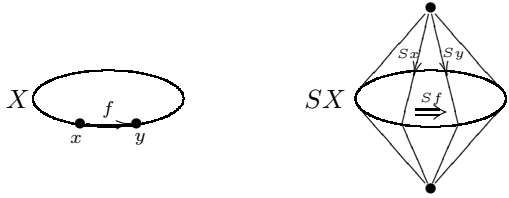
\bfig
\[
 \xy 
  (0,0)*\ellipse(10,3.75){-}; % space X.
 (-3.75,-3.35)*{\bullet}="S"+(-.5,-2)*{\scriptstyle x};
 (3.75,-3.35)*{\bullet}="S2"+(.5,-2)*{\scriptstyle y};
 {\ar@/_.05pc/^f "S";"S2"};
  (-12,0)*{X};
 \endxy
\qquad \qquad \qquad
 \xy
  (0,0)*\ellipse(10,3.75){-};
 (-3.75,-3.35)*{}="S";
 (3.75,-3.35)*{}="S2";
  (-9.75,0.9)*{}="LT";
  (-9.75,-0.9)*{}="LB";
  (9.75,0.9)*{}="RT";
  (9.75,-0.9)*{}="RB";
  (0,12)*{\bullet}="T";
  (0,-12)*{\bullet}="B";
  "T";"LT" **\dir{-};
  "T";"RT" **\dir{-};
  "B";"LB" **\dir{-};
  "RB";"B" **\dir{-};
  "T";"S" **\dir{-}; ?(.58)*\dir{>};
  "S";"B" **\dir{-};
  "T";"S2" **\dir{-}; ?(.58)*\dir{>};
  "S2";"B" **\dir{-};
  (-14,0)*{SX};
  (-2,-1.8)*{}="S";
 (2,-1.8)*{}="S2";
 {\ar@{=>}^<<<{\scriptscriptstyle Sf} "S";"S2"};
 (-3,6)*{\scriptscriptstyle Sx};
 (3,6)*{\scriptscriptstyle Sy};
 \endxy
\]
29.  Suspending $k$-morphisms in the fundamental $n$-groupoid
of a space $X$
\efig

For readers comfortable with homotopy theory we can make this a
bit more precise as follows.   Let $n\hot_k$ denote the category
of pointed $(k-1)$-connected homotopy $n+k$-types for $k \ge 1$,
or just homotopy $n$-types for $k = 0$.  Then we can think of
suspension as a functor $\tilde S \maps n\hot_k \to n\hot_{k+1}$,
where we take the reduced suspension for $k \ge 1$, while for $k
= 0$ we first adjoin a basepoint and then take the reduced
suspension.  Now suppose we have an equivalence $\Pi \maps
n\hot_k \to n\gpd_k$ to some category of weak
$n$-groupoids, with adjoint given by a functor $N \maps n\gpd \to
n\hot$, where $N$ is obtained by composing a
kind of `nerve' functor with `geometrical realization'.   One
expects that the $n$-categorical suspension functor of the
previous section yields one for weak $n$-groupoids, say $S \maps
n\gpd_k \to n\gpd_{k+1}$, where $n\gpd_k$ denotes the category of
$k$-tuply monoidal weak $n$-groupoids, and that $S$ is
equivalent to $N \tilde S \Pi$.  If this is so, the stabilization
result from homotopy theory should imply that $S$ is an
equivalence for $k \ge n+2$.    This argument also makes it clear
that the forgetful functor from $n\cat_k$ to $n\cat_{k-1}$ is
analogous to {\it looping} in homotopy theory \cite{KV}.

\section{Tangles}

Now let us venture an answer to the question:  {\it Of which
$n$-category are $n$-dimensional extended TQFTs representations?}
Recall from Section 2 that we expect this to be an $n$-category in
which the objects are 0-manifolds, the morphisms are 1-manifolds
with  boundary, the 2-morphisms are 2-manifolds with corners, and
so on.  These should all be oriented, but experience with
3-dimensional TQFTs has shown that it may be convenient to equip
them with some extra structure as well, a kind of `framing'.
One could attempt to make this more precise and give an explicit
description of the resulting $n$-category in terms of generators
and relations.  Instead, we try to isolate the crucial
algebraic properties of this $n$-category, and hypothesize that
it is in fact the universal $n$-category with these properties.
We have already seen that a key feature of the $n$-categories
appropriate for homotopy theory is the presence of {\it
inverses}.  Here we wish to argue that the corresponding feature
of the $n$-category of which extended TQFTs are representations
is the presence of {\it duals}.  Of course, both `inverses' and
`duals' may need to be understood in an appropriately weakened
sense.

We have already seen two levels of duality in the definition of a
unitary TQFT.  First, in $n\Cob$ each object $x$ is an oriented
$(n-1)$-manifold, and its dual $x^\ast$ is the same manifold with
its orientation reversed.  Second, each morphism $f \maps x \to
y$ is an oriented $n$-manifold with boundary, and its dual
$f^\dagger \maps y \to x$ is the same manifold with its
orientation reversed.   It is important to note that the dual
morphism $f^\dagger \maps y \to x$ is different from
the `adjoint' morphism $f^\ast \maps y^\ast \to x^\ast$ given by
the composite
\[
\begin{diagram}[xxxxxxx]
\node{y^\ast}
\arrow{e,t}{1 \tensor i_x}
\node{y^\ast \tensor x \tensor x^\ast}
\arrow{e,t}{1 \tensor f \tensor 1}
\node{y^\ast \tensor y \tensor x^\ast}
\arrow{e,t}{e_y \tensor 1}
\node{x^\ast.}
\end{diagram}
\]
The notion of adjoint morphism is derived from duality on
objects, but the notion of dual morphism is conceptually
independent.   Presumably this is just the tip of the iceberg:
in the $n$-categorical formulation of the theory of cobordisms
there should be $n+1$ distinct levels of duality, corresponding
to orientation reversal for $j$-manifolds with corners for all $0
\le j \le n$.

Indeed, while the details have only been worked out in special
cases, it seems that there should be a general
notion of a $k$-tuply monoidal $n$-category `with duals'.
Such a structure should have $n + 1$ duality operations, allowing us to take
duals of $j$-morphisms for all $0 \le j \le n$.  For the moment
let us denote all these duality operations with the symbol
${}^\star$.  For $j > 0$ the dual of a $j$-morphism $f \maps x \to
y$ is a $j$-morphism $f^\star \maps y \to x$, while the dual of a
0-morphism, or object, is simply another object.   For all $0 < j
< n$ there should be an associated unit and counit
\[   i_f
\maps 1_y \to ff^\star, \qquad
    e_f \maps f^\star f \to 1_x ,\]
satisfying the triangle identities, probably in a weakened sense.
Note that for $n$-morphisms there cannot be a  unit and counit
since there are no $(n+1)$-morphisms.   Also,  there can only be
a unit and counit for objects in the presence of a monoidal
structure, that is, if $k \ge 1$.   In this case, if we write the
monoidal structure as a tensor product, the unit and counit for
object take the usual forms as in eq.\  (\ref{unitcounit}). Other
properties we expect of duality are  $f^{\star \star} = f,$ and
$(fg)^\star = g^\star f^\star$ whenever the $j$-morphisms $f,g$
can be composed.

Since the precise axioms for a $k$-tuply monoidal $n$-category
`with duals' have only been formulated in certain semistrict
cases, let us go through them one by one.   We shall pay special
attention to the example of the {\it free semistrict $k$-tuply
monoidal $n$-category with duals on one object.}  We denote this
$n$-category by $C_{n,k}$ and the generating object by $x$.   We
hypothesize that the weak analog of this $n$-category
has a special significance for topological quantum field theory:
\vskip 1em\noindent
{\bf Tangle Hypothesis.}
The $n$-category of framed $n$-tangles in $n+k$ dimensions is
$(n+k)$-equivalent to the free weak $k$-tuply monoidal
$n$-category with duals on one object.
\vskip 1em\noindent

The sort of tangle studied in knot theory is what we would call a
1-tangle in 3 dimensions, shown in the $n = 1, k = 2$ entry of Figure
30 below.  By an $n$-tangle in $n+k$ dimensions we intend a
generalization of this concept, roughly speaking an $n$-manifold with
corners embedded in $[0,1]^{n+k}$ so that the codimension $j$ corners
of the manifold are mapped into the subset of $[0,1]^{n+k}$ for which
the last $j$ coordinates are either $0$ or $1$.  By a `framing' of an
$n$-tangle we mean a homotopy class of trivializations of the normal
bundle.  (Such a framing, together with the standard orientation on
$[0,1]^{n+k}$, determines an orientation on the submanifold.)

\bfig
\[ % 
 \xy
 (-75,60)*{k=0};
 (-75,30)*{k=1};
 (-75,0)*{k=2};
 (-75,-30)*{k=3};
 (-75,-60)*{k=4};
 (-40,75)*{n=0};
 (0,75)*{n=1};
 (40,75)*{n=2};
 (0,0)*{ %^^^^^^^^^^^^^^^^^^^^^^^^^^^^^^^^^^^^^^^^^^^^^^^^
   \xy 0;/r.20pc/: % n=1 k=2
   (25,25)*{};
 (-2,12)*{\bullet}="1"+(-2,1)*{ \scriptstyle x};
 (1,16)*{\bullet}="2"+(-2,2)*{\scriptstyle x^{\ast}};
 (10,14)*{\bullet}="3"+(1,3)*{\scriptstyle x};
 (1,-7)*{\bullet}="4"+(-1,3)*{\scriptstyle x};
 "3";"4" **\crv{} \POS?(.3)*{\hole}="J"; ?(0)*\dir{>};
 "1";"J" **\crv{(15,-5)}; ?(.15)*\dir{>};
 "2";"J" **\crv{} ?(.28)*\dir{<};
%*********************BUILDS BOX*****************************************************
 (-14,10)*{}="TL";
 (14,10)*{}="TR";
 (14,-10)*{}="BR";
 (-14,-10)*{}="BL";
 (-6,20)*{}="xTL";
 (22,20)*{}="xTR";
 (22,0)*{}="xBR";
 (-6,0)*{}="xBL";
     "TL";"TR" **\dir{-};
     "TR";"BR" **\dir{-};
     "BR";"BL" **\dir{-};
     "BL";"TL" **\dir{-};
     "xTL";"xTR" **\dir{-};
     "xTR";"xBR" **\dir{-};
     "xBR";"xBL" **\dir{.};
     "TL";"xTL" **\dir{-};
     "TR";"xTR" **\dir{-};
     "BL";"xBL" **\dir{.};
     "BR";"xBR" **\dir{-};
     "xTL";"xBL" **\dir{.};
 %*********************END BUILD BOX***************************************************
 \endxy
    };
    (0,30)*{
    %^^^^^^^^^^^^^^^^^^^^^^^^^^^^^^^^^^^^^^^^^^^^^^^^^^^^^^
     \xy 0;/r.20pc/: % n=1 k=1
 (-7,10)*{\bullet}="1"+(0,3)*{x};
 (-1,10)*{\bullet}="2"+(0,3)*{x^{\ast}};
 (6,10)*{\bullet}="3"+(0,3)*{x};
 (1,-10)*{\bullet}="4"+(0,-3)*{x};
  "1"; "2" **\crv{(-7,2) & (-1,2)}; ?(.2)*\dir{>}; ?(.9)*\dir{>};
  "3"; "4" **\crv{(9,4) & (-5,-1)}; ?(.5)*\dir{>};
 (-10,10)*{}; (10,10)*{} **\dir{-};
 (10,-10)*{}; (10,10)*{} **\dir{-};
 (-10,-10)*{}; (-10,10)*{} **\dir{-};
 (-10,-10)*{}; (10,-10)*{} **\dir{-};
 \endxy};
(40,0)*{
%^^^^^^^^^^^^^^^^^^^^^^^^^^^^^^^^^^^^^^^^^^^^^^^^^^^^^^
 \xy 0;/r.20pc/:% n=2 k=2
 (8.25,-1.25)*\ellipse(2,.65){-};
 (6,18)*{}="b";
   (1,16)*{}="a";
  \vunder~{(1,17.5)}{(6,18)}{(3.5,15.5)}{(6,16)};
  (1,17.5)*{}; (6,16)*{} **\crv{(-5,15)& (6,13) }; \POS?(.75)*{\hole}="J1";
  "J1";(6,18)*{} **\crv{(10,14)&(11,18)};
   (6,-2)*{}="b";
   (1,-4)*{}="a";
  \vunder~{(1,-2.5)}{(6,-2)}{(3.5,-4.5)}{(6,-4)};
  (1,-2.5)*{}; (6,-4)*{} **\crv{(-5,-5)& (6,-7) }; \POS?(.75)*{\hole}="J1";
  "J1";(6,-2)*{} **\crv{(10,-6)&(11,-2)};
  (-1,16)*{}="TL";
  (9,16)*{}="TR";
  (-1,-4)*{}="BL";
  (9.25,-3.5)*{}="BR";
  (18.5,-2.5)*{}="BRR";
  (14.5,-2.5)*{}="BRR2";
  "TL";"BL" **\crv{(-2,13) & (1,6)};
  "TR";"BRR" **\crv{(9,12) & (10,6)};
  "BR";"BRR2" **\crv{(7,12) & (12.5,0)};
%*********************BUILDS BOX*****************************************************
 (-14,10)*{}="TL";
 (14,10)*{}="TR";
 (14,-10)*{}="BR";
 (-14,-10)*{}="BL";
 (-6,20)*{}="xTL";
 (22,20)*{}="xTR";
 (22,0)*{}="xBR";
 (-6,0)*{}="xBL";
     "TL";"TR" **\dir{-};
     "TR";"BR" **\dir{-};
     "BR";"BL" **\dir{-};
     "BL";"TL" **\dir{-};
     "xTL";"xTR" **\dir{-};
     "xTR";"xBR" **\dir{-};
     "xBR";"xBL" **\dir{.};
     "TL";"xTL" **\dir{-};
     "TR";"xTR" **\dir{-};
     "BL";"xBL" **\dir{.};
     "BR";"xBR" **\dir{-};
     "xTL";"xBL" **\dir{.};
 %*********************END BUILD BOX***************************************************
 (5,22)*{\textbf{4d}};
 (25,25)*{};
 \endxy};
 (-40,0)*{ %^^^^^^^^^^^^^^^^^^^^^^^^^^^^^^^^^^^^^^^^
 \xy  0;/r.20pc/:% n=0 k=2
 (-5,5)*{\bullet}+(1,3)*{x};
 (5,4)*{\bullet}+(1,3)*{x^{\ast}};
 (1,-7)*{\bullet}+(1,3)*{x};
 (-10,10)*{}; (10,10)*{} **\dir{-};
 (10,-10)*{}; (10,10)*{} **\dir{-};
 (-10,-10)*{}; (-10,10)*{} **\dir{-};
 (-10,-10)*{}; (10,-10)*{} **\dir{-};
 \endxy
 };(0,-30)*{ %^^^^^^^^^^^^^^^^^^^^^^^^^^^^^^^^^^^^^^^^^^^^^^^^^
 \xy 0;/r.20pc/: % n=1 k=3
 (-4,14)*{\bullet}="1"+(,2)*{ \scriptstyle x};
 (0,14)*{\bullet}="2"+(,2)*{\scriptstyle x^{\ast}};
 (10,14)*{\bullet}="3"+(1,3)*{\scriptstyle x};
 (1,-7)*{\bullet}="4"+(-1,3)*{\scriptstyle x};
 "3";"4" **\crv{(10,4) & (-1,-6)}; ?(.6)*\dir{>};
%*********************BUILDS BOX*****************************************************
 (-14,10)*{}="TL";
 (14,10)*{}="TR";
 (14,-10)*{}="BR";
 (-14,-10)*{}="BL";
 (-6,20)*{}="xTL";
 (22,20)*{}="xTR";
 (22,0)*{}="xBR";
 (-6,0)*{}="xBL";
     "TL";"TR" **\dir{-}; \POS?(.35)*{\hole}="J"; \POS?(.5)*{\hole}="J1";
     "TR";"BR" **\dir{-};
     "BR";"BL" **\dir{-};
     "BL";"TL" **\dir{-};
     "xTL";"xTR" **\dir{-};
     "xTR";"xBR" **\dir{-};
     "xBR";"xBL" **\dir{.};
     "TL";"xTL" **\dir{-};
     "TR";"xTR" **\dir{-};
     "BL";"xBL" **\dir{.};
     "BR";"xBR" **\dir{-};
     "xTL";"xBL" **\dir{.};
 %*********************END BUILD BOX***************************************************
  "1";"J" **\crv{}?(.9)*\dir{>};
  "2";"J1" **\crv{} ?(.3)*\dir{<};
  "J";"J1" **\crv{(-4,5) & (0,5)};
  (5,22)*{\textbf{4d}};
  (25,25)*{};
 \endxy
 };
  (0,-60)*{%^^^^^^^^^^^^^^^^^^^^^^^^^^^^^^^^^^^^^^^^^^^^^^^^^^^^
    \xy 0;/r.20pc/: % n=1 k=4
 (-4,14)*{\bullet}="1"+(,2)*{ \scriptstyle x};
 (0,14)*{\bullet}="2"+(,2)*{\scriptstyle x^{\ast}};
 (10,14)*{\bullet}="3"+(1,3)*{\scriptstyle x};
 (1,-7)*{\bullet}="4"+(-1,3)*{\scriptstyle x};
 "3";"4" **\crv{(10,4) & (-1,-6)}; ?(.6)*\dir{>};
%*********************BUILDS BOX*****************************************************
 (-14,10)*{}="TL";
 (14,10)*{}="TR";
 (14,-10)*{}="BR";
 (-14,-10)*{}="BL";
 (-6,20)*{}="xTL";
 (22,20)*{}="xTR";
 (22,0)*{}="xBR";
 (-6,0)*{}="xBL";
     "TL";"TR" **\dir{-}; \POS?(.35)*{\hole}="J"; \POS?(.5)*{\hole}="J1";
     "TR";"BR" **\dir{-};
     "BR";"BL" **\dir{-};
     "BL";"TL" **\dir{-};
     "xTL";"xTR" **\dir{-};
     "xTR";"xBR" **\dir{-};
     "xBR";"xBL" **\dir{.};
     "TL";"xTL" **\dir{-};
     "TR";"xTR" **\dir{-};
     "BL";"xBL" **\dir{.};
     "BR";"xBR" **\dir{-};
     "xTL";"xBL" **\dir{.};
 %*********************END BUILD BOX***************************************************
  "1";"J" **\crv{}?(.9)*\dir{>};
  "2";"J1" **\crv{} ?(.3)*\dir{<};
  "J";"J1" **\crv{(-4,5) & (0,5)};
  (5,22)*{\textbf{5d}};
  (25,25)*{};
 \endxy
 };
 (-40,30)*{ %^^^^^^^^^^^^^^^^^^^^^^^^^^^^^^^^^^^^^^^^^^^^^^^^^^^^^^
 \xy   % n=0 k=1
 (0,0)*{\bullet}+(1,3)*{x^{\ast}};
 (5,0)*{\bullet}+(1,3)*{x};
 (-5,0)*{\bullet}+(1,3)*{x};
(-10,0)*{}; (10,0)*{} **\dir{-};
 \endxy
 };
 (-40,60)*{ %^^^^^^^^^^^^^^^^^^^^^^^^^^^^^^^^^^^^^^^^^^^^^^^^^^^^^^
% n=0 k=0
 \xy
 (0,0)*{\bullet}+(2,3)*{x^{\ast}};
 \endxy
 };
 (-40,-30)*{ %^^^^^^^^^^^^^^^^^^^^^^^^^^^^^^^^^^^^^^^^^^^^^^^^^^^^^^
   \xy 0;/r.20pc/: % n=0 k=3
 (-2.5,4)*{\bullet}+(1,2)*{\scriptstyle x};
 (5,7)*{ \bullet}+(3,1)*{\scriptstyle x^{\ast}};
 (6,0)*{\bullet}+(1,2)*{\scriptstyle x};
%*********************BUILDS BOX*****************************************************
 (-14,10)*{}="TL";
 (14,10)*{}="TR";
 (14,-10)*{}="BR";
 (-14,-10)*{}="BL";
 (-6,20)*{}="xTL";
 (22,20)*{}="xTR";
 (22,0)*{}="xBR";
 (-6,0)*{}="xBL";
     "TL";"TR" **\dir{-};
     "TR";"BR" **\dir{-};
     "BR";"BL" **\dir{-};
     "BL";"TL" **\dir{-};
     "xTL";"xTR" **\dir{-};
     "xTR";"xBR" **\dir{-};
     "xBR";"xBL" **\dir{.};
     "TL";"xTL" **\dir{-};
     "TR";"xTR" **\dir{-};
     "BL";"xBL" **\dir{.};
     "BR";"xBR" **\dir{-};
     "xTL";"xBL" **\dir{.};
 %*********************END BUILD BOX***************************************************
 (25,25)*{};
 \endxy
 };
 (-40,-60)*{
   \xy 0;/r.20pc/: % n=0 k=3
 (-2.5,4)*{\bullet}+(1,2)*{\scriptstyle x};
 (5,7)*{ \bullet}+(3,1)*{\scriptstyle x^{\ast}};
 (6,0)*{\bullet}+(1,2)*{\scriptstyle x};
%*********************BUILDS BOX*****************************************************
 (-14,10)*{}="TL";
 (14,10)*{}="TR";
 (14,-10)*{}="BR";
 (-14,-10)*{}="BL";
 (-6,20)*{}="xTL";
 (22,20)*{}="xTR";
 (22,0)*{}="xBR";
 (-6,0)*{}="xBL";
     "TL";"TR" **\dir{-};
     "TR";"BR" **\dir{-};
     "BR";"BL" **\dir{-};
     "BL";"TL" **\dir{-};
     "xTL";"xTR" **\dir{-};
     "xTR";"xBR" **\dir{-};
     "xBR";"xBL" **\dir{.};
     "TL";"xTL" **\dir{-};
     "TR";"xTR" **\dir{-};
     "BL";"xBL" **\dir{.};
     "BR";"xBR" **\dir{-};
     "xTL";"xBL" **\dir{.};
 %*********************END BUILD BOX***************************************************
 (25,25)*{};
(5,22)*{\textbf{4d}};
 \endxy
 };
 (0,60)*{
  \xy  0;/r.20pc/:  %k=o n=1
 (0,10)*{\bullet}="a"+(2.5,1)*{x};
 (0,-10)*{\bullet}="b"+(2.5,1)*{x};
 "a";"b" **\dir{-}; ?(.45)*\dir{>};
 \endxy
 };
 (40,60)*{
  \xy 0;/r.20pc/: % n=2 k=0
 (-10,10)*{\bullet}="TL"+(-1,3)*{x};
 (10,10)*{}="2"="TR"+(1,3)*{x};
 (10,-10)*{\bullet}="BR"+(1,-3)*{x};
 (-10,-10)*{\bullet}="BL"+(-1,-3)*{x};
 "TL";"TR" **\dir{-}; ?(.5)*\dir{>};
 "BL";"BR" **\dir{-}; ?(.5)*\dir{>};
 "TL";"BL" **\dir{-};
 "TR";"BR" **\dir{-};
 (0,3)*{}="a";
 (0,-3)*{}="a'";
 {\ar@{=>} "a";"a'"};
 \endxy
 };
 (40,30)*{
  \xy 0;/r.20pc/: % n=2 k=3
 (-9,10)*{\bullet}="1"+(-1,-3)*{x^{\ast}};
 (-1,10)*{\bullet}="2"+(-1,-3)*{x};
  (-9,-10)*{\bullet}="b1"+(-1,-3)*{x^{\ast}};
 (-1,-10)*{\bullet}="b2"+(-1,-3)*{x};
 (13,20)*{\bullet}="4"+(1,3)*{x};
 (7,20)*{\bullet}="3"+(1,3)*{x^{\ast}};
 (13,0)*{\bullet}="b4";
 (7,0)*{\bullet}="b3"+(1,3);
 "1";"2" **\crv{(-4,15) & (3,15)}; ?(.15)*\dir{<};
 "3";"4" **\crv{(3,15) & (10,15)}; ?(.85)*\dir{>};
  (-1,-5.25)*{}="M";
 "1";"b1"  **\dir{-};
 "2";"b2"  **\dir{-};
 "b1";"M"  **\dir{-}; ?(.76)*\dir{>};
 "b3";"M"  **\dir{.};
 "b2";"b4"  **\dir{-}; ?(.4)*\dir{<};
 "4";"b4"  **\dir{-};
 "3";"b3"  **\dir{.};
 (0,13.25)*{}="z1";
 (7,16.15)*{}="z2";
 "z1";"z2" **\crv{(2,-3) & (7,3)};
%*********************BUILDS BOX*****************************************************
 (-14,10)*{}="TL";
 (14,10)*{}="TR";
 (14,-10)*{}="BR";
 (-14,-10)*{}="BL";
 (-6,20)*{}="xTL";
 (22,20)*{}="xTR";
 (22,0)*{}="xBR";
 (-6,0)*{}="xBL";
     "TL";"TR" **\dir{-};
     "TR";"BR" **\dir{-};
     "BR";"BL" **\dir{-};
     "BL";"TL" **\dir{-};
     "xTL";"xTR" **\dir{-};
     "xTR";"xBR" **\dir{-};
     "xBR";"xBL" **\dir{.};
     "TL";"xTL" **\dir{-};
     "TR";"xTR" **\dir{-};
     "BL";"xBL" **\dir{.};
     "BR";"xBR" **\dir{-};
     "xTL";"xBL" **\dir{.};
 %*********************END BUILD BOX***************************************************
 (25,25)*{};
 \endxy
 };
 (40,-30)*{
  \xy 0;/r.20pc/:% n=2 k=3
 (-1.5,-2)*\ellipse(3,1){-};
 (3.5,-2)*\ellipse(3,1){-};
 (6,18)*{}="b";
   (1,16)*{}="a";
  \vunder~{(1,17.5)}{(6,18)}{(3.5,15.5)}{(6,16)};
  (1,17.5)*{}; (6,16)*{} **\crv{(-5,15)& (6,13) }; \POS?(.75)*{\hole}="J1";
  "J1";(6,18)*{} **\crv{(10,14)&(11,18)};
  (-1,16)*{}="TL";
  (9,16)*{}="TR";
  (-6,-4)*{}="BLL";
  (0,-4)*{}="BL";
  (10,-4)*{}="BRR";
  (4,-4)*{}="BR";
  (3,5)*{}="C";
  (4.25,6.25)*{}="C2";
  (4.25,15)*{}="C3";
    "TL";"BLL" **\crv{(-2,13) & (-1,6)};
    "TR";"BRR" **\crv{(9,13) & (10,6)};
    "C";"BL" **\crv{};
    "C";"BR" **\crv{(6,10)};
     "C3";"C2" **\dir{.};
%*********************BUILDS BOX*****************************************************
 (-14,10)*{}="TL";
 (14,10)*{}="TR";
 (14,-10)*{}="BR";
 (-14,-10)*{}="BL";
 (-6,20)*{}="xTL";
 (22,20)*{}="xTR";
 (22,0)*{}="xBR";
 (-6,0)*{}="xBL";
     "TL";"TR" **\dir{-};
     "TR";"BR" **\dir{-};
     "BR";"BL" **\dir{-};
     "BL";"TL" **\dir{-};
     "xTL";"xTR" **\dir{-};
     "xTR";"xBR" **\dir{-};
     "xBR";"xBL" **\dir{.};
     "TL";"xTL" **\dir{-};
     "TR";"xTR" **\dir{-};
     "BL";"xBL" **\dir{.};
     "BR";"xBR" **\dir{-};
     "xTL";"xBL" **\dir{.};
 %*********************END BUILD BOX***************************************************
 (5,22)*{\textbf{5d}};
 (25,25)*{};
 \endxy
 };
 (40,-60)*{
\xy 0;/r.20pc/:% n=2 k=4
 (-1.5,-2)*\ellipse(3,1){-};
 (3.5,-2)*\ellipse(3,1){-};
 (6,18)*{}="b";
   (1,16)*{}="a";
  \vunder~{(1,17.5)}{(6,18)}{(3.5,15.5)}{(6,16)};
  (1,17.5)*{}; (6,16)*{} **\crv{(-5,15)& (6,13) }; \POS?(.75)*{\hole}="J1";
  "J1";(6,18)*{} **\crv{(10,14)&(11,18)};
  (-1,16)*{}="TL";
  (9,16)*{}="TR";
  (-6,-4)*{}="BLL";
  (0,-4)*{}="BL";
  (10,-4)*{}="BRR";
  (4,-4)*{}="BR";
  (3,5)*{}="C";
  (4.25,6.25)*{}="C2";
  (4.25,15)*{}="C3";
    "TL";"BLL" **\crv{(-2,13) & (-1,6)};
    "TR";"BRR" **\crv{(9,13) & (10,6)};
    "C";"BL" **\crv{};
    "C";"BR" **\crv{(6,10)};
     "C3";"C2" **\dir{.};
%*********************BUILDS BOX*****************************************************
 (-14,10)*{}="TL";
 (14,10)*{}="TR";
 (14,-10)*{}="BR";
 (-14,-10)*{}="BL";
 (-6,20)*{}="xTL";
 (22,20)*{}="xTR";
 (22,0)*{}="xBR";
 (-6,0)*{}="xBL";
     "TL";"TR" **\dir{-};
     "TR";"BR" **\dir{-};
     "BR";"BL" **\dir{-};
     "BL";"TL" **\dir{-};
     "xTL";"xTR" **\dir{-};
     "xTR";"xBR" **\dir{-};
     "xBR";"xBL" **\dir{.};
     "TL";"xTL" **\dir{-};
     "TR";"xTR" **\dir{-};
     "BL";"xBL" **\dir{.};
     "BR";"xBR" **\dir{-};
     "xTL";"xBL" **\dir{.};
 %*********************END BUILD BOX***************************************************
 (5,22)*{\textbf{6d}};
 (25,25)*{};
 \endxy
 }:
 \endxy % END END END END END END END END END END
\]
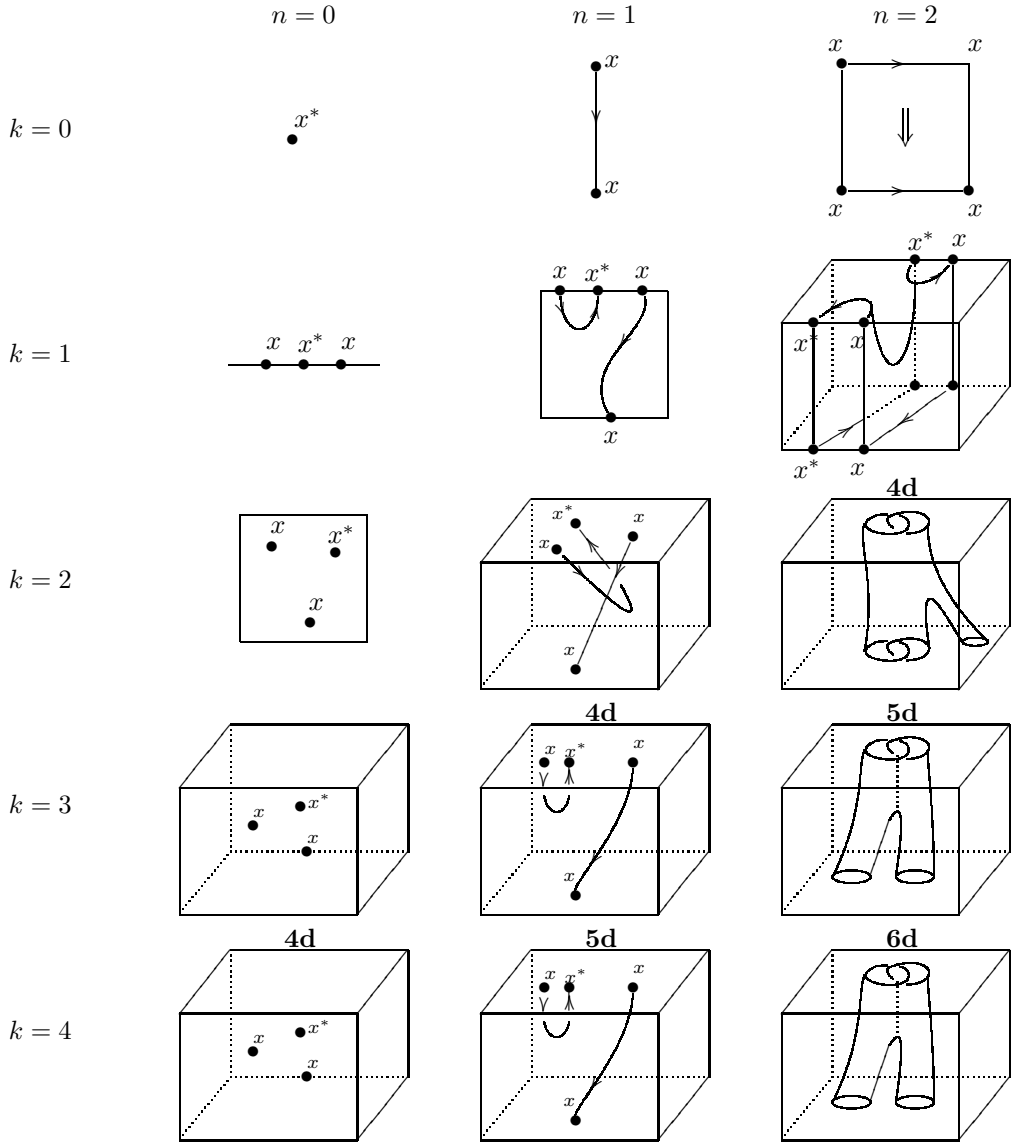
30.  Examples of $n$-morphisms in $C_{n,k}$, drawn as
$n$-tangles in $n+k$ dimensions
\efig

Implicit in the tangle hypothesis is that there should be a weak
$n$-category whose $n$-morphisms are suitably defined isotopy classes
of $n$-tangles in $n+k$ dimensions.  The precise definition of
`$n$-tangles' and `isotopy classes' for this purpose is only
well-understood for $n = 1$, and only beginning to be understood for
$n = 2$.  In what sense should isotopy classes of $n$-tangles in $n+k$
dimensions be the $n$-morphisms in an $n$-category?  Roughly, a
$j$-morphism $f \maps x \to y$ in this $n$-category should be a
certain equivalence class of $j$-tangles in $[0,1]^{j+k}$, going from
the equivalence class $x$ of $(j-1)$-tangles in $[0,1]^{j+k-1} \times
\{0\}$ to the equivalence class $y$ of $(j-1)$-tangles in
$[0,1]^{j+k-1}\times \{1\}$.  The duality operation on $j$-tangles
corresponds to reflection of $[0,1]^{j+k}$ about the last coordinate
axis.

Clearly much remains to be made precise here. Rather than
continuing to speak in generalities, let us describe what is
known so far in various special cases.   This should also be the
best way of illustrating the significance of the $k$-tuply monoidal
structure.

We begin with $n = 0$ column of Figure 30.  Here there is just
one level of duality to consider, which we denote by $\ast$.  For $k
= 0$, a $k$-tuply monoidal $n$-category $C$ is just a set.  A set
`with duals' is simply one equipped with an involution, that is,
a function $\ast \maps C \to C$ with $x^{\ast\ast} = x$ for all
$x$.  Thus $C_{0,0}$ is the free set with involution on one
object $x$, namely the two-element set $\{x,x^\ast\}$.    Now for
$k > 0$, a framing of a 0-dimensional submanifold of $[0,1]^k$ is
equivalent to an orientation, so in this degenerate case we
somewhat artificially identify a framing with an orientation.
The cobordism hypothesis then states that $C_{0,0}$ describes the
set of isotopy classes of oriented 0-tangles in 0 dimensions ---
i.e., the positively and negatively oriented point!

Proceeding down the column to $k = 1$, a $k$-tuply monoidal
$n$-category $C$ is just a monoid.  A monoid `with duals' is one
equipped with an involution, which in this context means a
function $\ast \maps C \to C$ with $x^{\ast\ast} = x$ and $(xy)^\ast =
y^\ast x^\ast$.   Thus $C_{0,1}$ is the free monoid with
involution on one object $x$.  Elements of $C_{0,1}$ are thus
formal (noncommuting) products  of the elements $x$ and $x^\ast$.
These correspond to isotopy classes of oriented
$0$-tangles in $[0,1]$, or simply strings of positively and
negatively oriented points, as shown in Figure 30.

Continuing down to $k = 2$, a $k$-tuply monoidal $n$-category is
a commutative monoid, and again `having duals' simply means being
equipped with an involution.   Thus $C_{0,2}$ is the free
commutative monoid on one object $x$.  Any element of $C_{0,2}$
is thus of the form $x^n (x^\ast)^m$.   Similarly, an isotopy
class of oriented $0$-tangles in 2 dimensions is just a
collection of $n$ positively oriented and $m$ negatively oriented
points.  The same holds for all dimensions $k \ge 2$, essentially
because there are enough dimensions to freely move the points
about in a manner corresponding to the Eckmann-Hilton argument.

In the $n = 1$ column we have various kinds of category with
extra structure.   Here there are two levels of duality ---
duality for objects and duals for morphisms --- so to avoid
confusion we denote the dual of an object $x$ by $x^\ast$, and
the dual of a morphism $f \maps x \to y$ by $f^\dagger \maps y
\to x$.  For $k = 0$ we just have categories, and by a category
`with duals' we mean one equipped with operations $\ast$ and
$\dagger$ as just described, such that $\ast^2 = 1$ on objects,
$\dagger^2 = 1$ on morphisms, and $(fg)^\dagger = g^\dagger
f^\dagger$.  With this definition, $C_{1,0}$, the free category
with duals on one object $x$, is rather dull: it has objects $x$
and $x^\ast$ and only identity morphisms.   If we artificially
take a framing in this degenerate case to mean an orientation,
the tangle hypothesis states that the morphisms in this category
correspond to isotopy classes of oriented  1-tangles in 1
dimension.  There are indeed only two of these, one of which is
shown in Figure 30.

Moving down the $n = 1$ column to $k = 1$, we have monoidal
categories.  This is the first point at which there is room for
the unit and counit.  By a monoidal category `with duals' we thus
mean a category with duals in the above sense which is also
monoidal, for which $(x \tensor y)^\ast = y^\ast \tensor x^\ast$
and $(f \tensor g)^\dagger = f^\dagger \tensor g^\dagger$,  and
equipped with a unit and counit for the $\ast$ duality, natural
morphisms satisfying the triangle identity.   We also require
\[   e_x^\dagger = i_{x^\ast}.\]
This expresses a relation between the two levels of duality, but
we should emphasize that the relations between different levels
of $n$-categorical duality are very poorly understood, and our
treatment here is provisional.  In any event, with this
definition the morphisms in $C_{1,1}$ describe isotopy classes of
framed (or equivalently, oriented) 1-tangles in 2 dimensions
\cite{FY}. Composition of morphisms in this category corresponds
to  vertical juxtaposition of 1-tangles, while the tensor product
corresponds to horizontal juxtaposition.

Continuing down to $k = 2$, we have braided monoidal categories.
By a braided monoidal category `with duals' we mean a monoidal
category with duals which is also braided.  We also require
that the braiding and the `balancing' shown in Figure 31 be
unitary, where a morphism $f$ is said to be unitary if
$f f^\dagger = f^\dagger f = 1$.

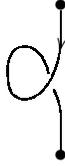
\bfig 
\[
 \xy 
 (0,10)*{\bullet}="T";
 (0,-10)*{\bullet}="B";
 (0,5)*{}="T'";
 (0,-5)*{}="B'";
 "T";"T'" **\dir{-};
 "B";"B'" **\dir{-};
 (-3,0)*{}="MB";
 (-7,0)*{}="LB";
    "T'";"LB" **\crv{(-1,-4) & (-7,-4)}; \POS?(.25)*{\hole}="2z"; ?(.05)*\dir{>};
    "LB"; "2z" **\crv{(-8,6) & (-2,6)};
    "2z"; "B'"  **\crv{(0,-3)};
 \endxy
\]
31. The balancing in a braided monoidal category with duals
\efig

Turaev and Yetter \cite{Turaev,Yetter} have shown that the morphisms
in $C_{1,2}$ correspond to isotopy classes of framed 1-tangles in 3
dimensions.  Here a couple of remarks are in order.  First, in this
dimension, our sort of framing is equivalent to an orientation
together with what is often called a framing in knot theory.  Second,
it is common to study 1-tangles in 3 dimensions using categories in
which the existence of the balancing is a separate postulate
\cite{CP}.  In our approach it arises automatically.  In fact, this
idea is already implicit in the work of Fr\"ohlich and Kerler
\cite{FK}.

By a symmetric monoidal category `with duals' we mean just a
braided monoidal category with duals which is also symmetric.
The morphisms in $C_{1,3}$ correspond to isotopy classes of
framed 1-tangles in $4$ dimensions, because $C_{1,3}$ is obtained
from $C_{1,2}$ by  adding the extra relations $B_{u,v} =
B_{v,u}^{-1}$ for all $u,v$, corresponding to the fact that in 4
dimensions there is room to unlink all links.  The same is true
in all higher dimensions --- as one would expect from the
stabilization hypothesis.  More generally, transversality results
from differential topology imply that for $k \ge n+2$, all
embeddings of compact $n$-manifolds in $\R^{n+k}$ are isotopic
\cite{Hirsch}.  Given the tangle hypothesis, this is a powerful
piece of evidence for the stabilization hypothesis.

In the $n = 2$ column we find the situation less well understood
but still very promising.  The precise definition of a $k$-tuply
monoidal 2-category with duals has not yet been systematically
worked out.   Instead, work so far has focused on the relation
between braided monoidal 2-categories with duals and 2-tangles in
4 dimensions.   Carter, Saito and others have worked out a
description of such 2-tangles as movies in which each frame is a
1-tangle in 3 dimensions, giving explicit `movie moves' which go
between any two movies representing isotopic 2-tangles \cite{CS}.
Fischer \cite{Fischer} has used this information to describe a
2-category of 2-tangles in 4 dimensions, and came close to
proving the tangle hypothesis in this case.   There are a number
of loose ends, however, and recently Kharlamov and Turaev
\cite{KT} have done some careful work on this subject,
particularly concerning the equivalence relation on 1-tangles in
3 dimensions needed to defining 1-morphisms as equivalence
classes.

Rather than reviewing this work in detail, let us simply touch upon
what may at first seem the most surprising feature, namely, how
the various levels of duality interact to yield the 2-morphisms
corresponding to the movies in Figure 6: the birth of a
circle, death of a circle, and saddle point.   This after all, is a
crucial part of the tangle conjecture: that duality in $n$-categories
naturally yields the correct handle attachments and handle
cancellations in the Morse-theoretic description of $n$-tangles.

We expect the movies in Figure 6 to occur whenever $k \ge 1$, so
for simplicity consider the case $k = 1$.
In a monoidal 2-category with duals
there should be 3 levels of duality: each object $x$ should have
a dual $x^\ast$,  each morphism $f \maps x \to y$ should have a
dual $f^\dagger \maps y \to x$, and each 2-morphism $\alpha \maps
f \doublearrow g$ should have a dual $\hat \alpha \maps g
\doublearrow f$.   The dualities for objects and morphisms should
have associated units and counits, but not the duality for
2-morphisms.  Thus for any object $x$ there are morphisms
\[      i_x \maps 1 \to x \tensor x^\ast, \qquad
        e_x \maps x^\ast \tensor x \to 1 ,\]
and for any morphism $f \maps x \to y$ there are 2-morphisms
\[   \iota_f \maps 1_y \doublearrow ff^\dagger, \qquad
     \epsilon_f \maps f^\dagger f \doublearrow 1_x .\]
The unit and counit for morphisms should satisfy the triangle
identities, while the unit and counit for objects should probably satisfy
them only weakly, i.e., up to the 2-isomorphisms shown as
2-tangles in 3 dimensions in Figure 32.

\bfig 
\[
 \xy 
 (-8,-15)*{\bullet}="BB"+(-1,-3)*{x};
 (-8,0)*{\bullet}="BT"+(-1,3)*{x};
 (10,-5)*{\bullet}="TB"+(1,-3)*{x};
 (10,10)*{\bullet}="TT"+(1,3)*{x};
 (0,-3)*{}="F";
 "BT";"TT" **\dir{-}; ?(.45)*\dir{<};
 "TB";"F" **\crv{(4,-3) & (-5, -25)} \POS?(.92)*{\hole}="J1"; \POS?(.78)*{\hole}="J2";
 "BB";"BT" **\dir{-};
 "TT";"TB" **\dir{-};
 "BB";"J1" **\crv{(-5,-13) & (-4,-5)}; ?(.25)*\dir{<};
  "J2";"J1" **\crv{~*=<2pt>{.} (3,-7)};
 \endxy
\qquad \qquad \qquad
 \xy
 (-8,-15)*{\bullet}="BB"+(-1,-3)*{x^{\ast}};
 (-8,0)*{\bullet}="BT"+(-1,3)*{x^{\ast}};
 (10,-5)*{\bullet}="TB"+(1,-3)*{x^{\ast}};
 (10,10)*{\bullet}="TT"+(1,3)*{x^{\ast}};
 (0,-3)*{}="F";
  "BB";"BT" **\dir{-};
 "TT";"TB" **\dir{-};
 "BT";"TT" **\dir{-}; ?(.55)*\dir{>};
 "BB";"F" **\crv{(-5,-7) & (14, -20)} \POS?(.92)*{\hole}="J2";
 \POS?(.75)*{\hole}="J1";
 "TB";"J2" **\crv{}; ?(.5)*\dir{<};
  "J1";"J2" **\crv~pC{~*=<2pt>{.} (-7,-8)};
 \endxy
\]
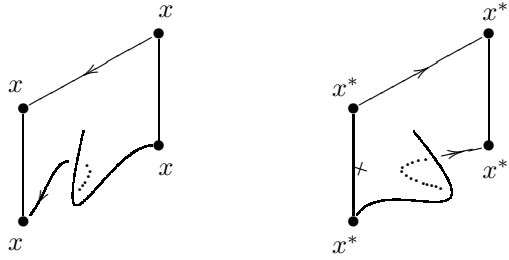
32.  The 2-isomorphisms $1_x \doublearrow (1_x \tensor e_x)(i_x
\tensor 1_x)$ and $1_{x^\ast} \doublearrow (e_x \tensor
1_{x^\ast})(1_{x^\ast} \tensor i_x)$ \break as 2-tangles in 3 dimensions
\efig

We also expect rules such as
\[    (x \tensor y)^\ast = y^\ast \tensor x^\ast, \qquad
       (fg)^\dagger = g^\dagger f^\dagger, \qquad
       (\alpha\beta)\hat{} = \hat\beta \hat\alpha ,\]
\[    x^{\ast\ast} = x, \qquad f^{\dagger\dagger} = f, \qquad
\hat{\hat\alpha} = \alpha, \]
and
\[   e_x^\dagger = i_{x^\ast},\qquad \hat \eps_f = \iota_{f^\dagger}.\]

Then, as shown in Figure 33, the birth of a clockwise oriented
circle corresponds to the unit of the counit of $x$, that is,
$\iota_{e_x}  \maps 1_1 \doublearrow e_x e_x^\dagger $.
Similarly, the death of a clockwise oriented circle corresponds
to the counit  of the unit of $x$, $\eps_{i_x} \maps i_x^\dagger
i_x \doublearrow 1_1 $. One sort of saddle point is given by the
unit of the unit of $x$, $\iota_{i_x} \maps 1_{x \tensor x^\ast}
\doublearrow i_x i_x^\dagger ,$ while another is given by the
counit of the counit of $x$,  $\eps_{e_x} \maps  e_x^\dagger e_x
\doublearrow 1_{x^\ast \tensor x}$.  Differently oriented
versions of the 1-manifolds with boundary appearing here can be
obtained by replacing $x$ with $x^\ast$ above.

\bfig
\[ 
 \vcenter{\xy
    (0,-4)*\ellipse(5,2){.};
    (0,-4)*\ellipse(5,2)__,=:a(-180){-};
    (-5,-8)*{}="TL";
    (5,-8)*{}="TR";
     "TL";"TR" **\crv{(-5,4) & (5,4) };
     (0,-10)*{\scriptstyle <};
 \endxy}
\qquad
 \vcenter{\xy
    (0,4)*\ellipse(5,2){-};
    (0,18)*{};
    (-5,8)*{}="TL";
    (5,8)*{}="TR";
     "TL";"TR" **\crv{(-5,-4) & (5,-4) };
     (0,6)*{\scriptstyle <};
 \endxy}
\qquad
 \xy 0;/r.18pc/:
  (20,2)*{\bullet}="RU"+(1,3)*{x^{\ast}};
  (16,-3)*{\bullet}="RD"+(2.5,2)*{x^{\ast}};
  (-16,2)*{\bullet}="LU"+(-1,2)*{x};
  (-20,-3)*{\bullet}="LD"+(-1,2)*{x};
  "RU";"RD" **\crv{(4,2) & (4,-1)}; ?(.1)*\dir{<}; ?(.85)*\dir{<};
  "LD";"LU" **\crv{(-4,-2) & (-4,1)}; ?(.08)*\dir{<}; ?(.85)*\dir{<};
    (7.5,0)*{}="x1";
    (-7.5,0)*{}="x2";
     "x1"; "x2" **\crv{(7,-10) & (-7,-10)};
  (16,-20)*{\bullet}="RDD"+(2.5,-1)*{x^{\ast}};
  (-20,-20)*{\bullet}="LDD"+(-1,-3)*{x};
   (20,-12.5)*{\bullet}="RUD"+(3.5,1)*{x^{\ast}};
   (-16,-15)*{\bullet}="LUD";
   (-16,-2.5)*{}="A";
   (16.1,-14.9)*{}="B";
        "RD"; "RDD" **\dir{-};
        "LD"; "LDD" **\dir{-};
        "A"; "LUD" **\dir{.};
        "RDD"; "LDD" **\crv{(0,-17)}; ?(.57)*\dir{>};
        "RU"; "RUD" **\dir{-};
        "LU"; "A" **\dir{-};
        "B"; "RUD" **\crv{(18,-14.15)};
        "B"; "LUD" **\crv{~*=<4pt>{.}(0,-18)}; ?(.4)*\dir{<};

 \endxy
 \qquad
 \xy 0;/r.18pc/:
 (-10,-12)*{\bullet}="A"+(-3,1)*{x};
 (-20,-6)*{\bullet}="B"+(-3,1)*{x^{\ast}};
  (10,10)*{\bullet}="a"+(1,3)*{x^{\ast}};
 (20,6)*{\bullet}="b"+(3,1)*{x};
  (-10,-24)*{\bullet}="A1"+(-3,1)*{x};
 (-20,-20)*{\bullet}="B1"+(-3,1)*{x^{\ast}};
  (10,2)*{}="a1";
  (10,-4)*{}="a11";
 (20,-6)*{\bullet}="b1"+(3,-1)*{x};
 "A";"B" **\crv{(0,0) & (-10,0)}; ?(.85)*\dir{<}; ?(.05)*\dir{<};
 "a";"b" **\crv{(0,0) & (10,0)};?(.85)*\dir{<}; ?(.05)*\dir{<};
 (-6,-3)*{}="X";
 (6,3)*{}="x";
   (-10,-14.5)*{}="P";
  "X";"x" **\crv{(1,-16)};
  "A";"A1" **\dir{-};
  "B";"B1" **\dir{-};
  "b";"b1" **\dir{-};
  "a";"a1" **\dir{-};
  "a11";"a1" **\dir{.};
  "P";"a11" **\dir{.};
  "B1";"P" **\dir{-}; ?(.55)*\dir{>};
  "A1";"b1" **\dir{-}; ?(.45)*\dir{<};
 \endxy
 \]
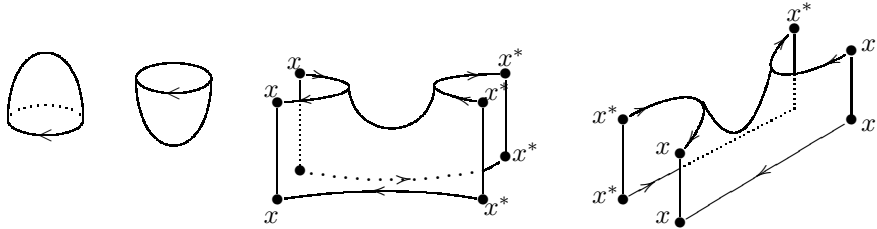
33.  The 2-morphisms $\iota_{e_x}, \eps_{i_x}, \iota_{i_x},$ and
$\eps_{e_x}$ as 2-tangles in 3 dimensions
\efig
\noindent
The reader may enjoy studying how the triangle identities for
$\iota$ and $\eps$ translate into handle cancellations. For
example, the following triangle identity, where $\tensor$ denotes
horizontal composition:
\[
\begin{diagram} [i_x i_x^\dagger i_x]
\node{i_x} \arrow[2]{e,t}{1_{i_x}} \arrow{se,b}{\iota_{i_x} \tensor 1_{i_x}}
\node[2]{i_x}   \\
\node[2]{i_x i_x^\dagger i_x} \arrow{ne,r}{1_{i_x} \tensor \eps_{i_x}}
\end{diagram}
\]
gives the isotopy between 2-tangles in 3 dimensions shown in
Figure 34.

\bfig
\[ % FIG.32 A 2-categorical triangle identity as 2-tangle isotropy
 \xy 0;/r.18pc/:
 (-10,-2)*{\bullet}="A"+(-3,1)*{x};
 (-20,4)*{\bullet}="B"+(-3,1)*{x^{\ast}};
  (-10,-24)*{\bullet}="A1"+(-3,1)*{x};
 (-20,-20)*{\bullet}="B1"+(-3,1)*{x^{\ast}};
 "A";"B" **\crv{(0,10) & (-10,10)}; ?(.85)*\dir{<}; ?(.05)*\dir{<};
  "A1";"B1" **\crv{(0,-12) & (-10,-12)}; ?(.85)*\dir{<}; ?(.05)*\dir{<};
  "A";"A1" **\dir{-};
  "B";"B1" **\dir{-};
  (-6,6.75)*{}="x1";
  (-6,-15)*{}="x2";
  "x1";"x2" **\dir{-};
 \endxy
\qquad =\qquad
 \xy 0;/r.18pc/:
 (-10,-2)*{\bullet}="A"+(-3,1)*{x};
 (-20,4)*{\bullet}="B"+(-3,1)*{x^{\ast}};
  (-10,-24)*{\bullet}="A1"+(-3,1)*{x};
 (-20,-20)*{\bullet}="B1"+(-3,1)*{x^{\ast}};
 "A";"B" **\crv{(0,10) & (-10,10)}; ?(.85)*\dir{<}; ?(.05)*\dir{<};
  "A1";"B1" **\crv{(0,-12) & (-10,-12)}; ?(.85)*\dir{<}; ?(.05)*\dir{<};
  "A";"A1" **\dir{-};
  "B";"B1" **\dir{-};
  (-6,6.75)*{}="x1";
  (-6,-15)*{}="x2";
  "x1";"x2" **\crv{(-8,-6) & (8,10) & (7,-9) & (-9,-11)};
 \endxy
 \]
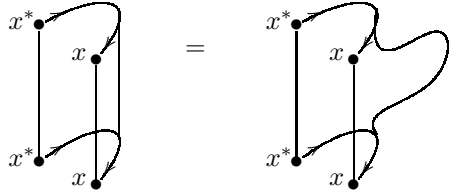
34.  A 2-categorical triangle identity as 2-tangle isotopy
\efig

Proceeding down the $n = 2$ column in Figure 30, we
observe the following features.  Taking $k = 2$, a 2-morphism
$\alpha\maps f \to g$ in $C_{2,2}$  should correspond to a
2-tangle in 4 dimensions, going from an equivalence class $f$ of
1-tangles in 3 dimensions to another equivalence class $g$.  The
braiding phenomena arise from the fact that $C_{2,2}$ is
a braided monoidal 2-category \cite{CS,Fischer}.  The case $k = 3$
corresponds to 2-tangles in 5 dimensions.   In Figure 30 we have
shown such a 2-tangle going from a seemingly linked pair of
circles to an unlinked pair.   The point is that, thinking of the 5th
dimension as time, any link in 4-dimensional space can be unlinked as
time passes.  In fact, an over-crossing can be isotoped to an
under-crossing while pushing one strand either a little bit
`up' into the 4th dimension or a little bit `down'.
Algebraically, since $C_{2,3}$ is weakly involutory, we expect
these two distinct isotopies to correspond to
$I_{x,y} \maps B_{x,y} \doublearrow
B_{y,x}^{-1}$ and $(I_{x,y}^{-1})_{hor}^{-1} \maps B_{x,y}
\doublearrow B_{y,x}^{-1}$.  The case $k = 4$
corresponds to 2-tangles in 6 dimensions.  In this dimension and
higher dimensions, the two isotopies from an over-crossing
to an under-crossing are themselves isotopic.  This should correspond
to the fact that in a strongly involutory monoidal 2-category,
$I_{x,y} = (I_{x,y}^{-1})_{hor}^{-1}$.

Now, taken together, the stabilization hypothesis and tangle
hypothesis suggest a proposal for the
$n$-category of which $n$-dimensional extended TQFTs are representations:
\vskip 1em \noindent
{\bf Extended TQFT Hypothesis, Part I.}
The $n$-category of which $n$-dimensional extended TQFTs are
representations is
the free stable weak $n$-category with duals on one object.
\vskip 1em\noindent

The idea here is, first, that by the stabilization hypothesis
$C_{n,k}$ should stabilize for $k \ge n+2$, yielding what we call the
`free stable weak $n$-category with duals on one object',
$C_{n,\infty}$.  Topologically this suggests that when considering the
category whose objects are framed 0-manifolds, whose morphisms
are framed 1-manifolds with boundary, whose 2-morphisms are
framed 2-manifolds with corners, and so on up to $n$, we are
free to think of all these objects as embedded in $[0,1]^{n+k}$
where $k \ge n+2$.  Thus in particular the $n$-morphisms can
be thought of as isotopy classes of framed $n$-tangles in $n+k$
dimensions, where $k \ge n+2$.

In the next section we turn to the implications of this
hypothesis for topological quantum field theory; here we briefly
summarize one more sophisticated piece of evidence for it.
This comes from the connection between framed cobordism theory
and  stable homotopy theory \cite{Adams,Stong}.  By the relation
between the $n$-categorical and homotopy-theoretic notions of
suspension, one expects the homotopy $n$-type of  $\Omega^k S^k$
to correspond to a very special weak $n$-groupoid, namely the
free $k$-tuply monoidal weak $n$-groupoid on one object,
$G_{n,k}$.  Now duals are simply a weakened form of inverses
\cite{KVinfinity}, so by the universal property of $C_{n,k}$ one
expects there to be a weak $n$-functor $T \maps C_{n,k} \to
G_{n,k}$ that turns duals into inverses.   Topologically
speaking, $T$ should be given by the Thom-Poyntragin
construction.  Indeed, this construction is implicit in
Figure 26, where a tangle is used to describe a homotopy.
Suppose, for example, that $\alpha = \beta$ is the generating
object $x$ of $C_{1,2}$.  Then the tangle in Figure 26 is the morphism
$B_{x,x}$ in $C_{1,2}$.   On the other hand, $Tx$ is the object
in $G_{1,2}$, or point in $\Omega^2 S^2$, corresponding to
identity map from $S^2$ to itself.  Thus $TB_{x,x}$ is a morphism
in $G_{1,2}$ corresponding to a nontrivial homotopy from $Tx
\tensor Tx$ to itself.

More generally, define a `$j$-loop' in any monoidal $n$-category
to be an $j$-morphism from $1_{j-1}$ to $1_{j-1}$, where we
define the object $1_0$ to be the unit for the monoidal structure
and, recursively, $1_{i+1} = 1_{1_i}$.  Now $T$ should map
$j$-loops in $C_{n,k}$ to $j$-loops in $G_{n,k}$.  In particular,
an $n$-loop in $C_{n,k}$ is just an isotopy class of compact
framed $j$-manifolds embedded in $[0,1]^{n+k}$ (or equivalently
$S^{n+k}$).   On the other hand, an $n$-loop in $G_{n,k}$ is an
an element of $\pi_{n+k}(S^k)$.  There is indeed a map from the
former to the latter, given by the Thom-Pontryagin construction.
Now for $k \ge n+2$, $G_{n,k}$ should stabilize to a weak
$n$-groupoid $G_{n,\infty}$ representing the homotopy $n$-type of
the infinite loop space $\Omega^\infty S^\infty$, so we expect to
obtain a weak $n$-functor $T_\infty \maps C_{n,\infty} \to
G_{n,\infty}$.   By the universal property of $G_{n,\infty}$,
this $n$-groupoid should simply be the result of adjoining
formal inverses to all $j$-morphisms in $C_{n,\infty}$.   One
expects from this that the $i$th framed cobordism group is
isomorphic to $\pi_i(\Omega^\infty S^\infty)$, that is, the $i$th
stable homotopy group of spheres.  This is indeed the case!

\section{Extended TQFTs}

One can think of $n$-category theory as providing a natural
hierarchy of generalizations of set theory.   The
basic idea is that the mathematics of sets, regarded
as the study of the category $\Set$, leads us
to consider general categories.  Regarding this as the
study of the 2-category $\Cat$ we are then lead to consider
general 2-categories, and so on.  In general, the category
$n\Cat$ of small strict $n$-categories is a strict
$(n+1)$-category, and we expect something of a similar but more
sophisticated sort to hold for weak $n$-categories.
At each level of this hierarchy one can do abstract algebra,
which at at the $n$th level is intimately tied to $n$-dimensional
topology.   To describe an extended TQFT as a representation of
an $n$-category, we must also develop analogs of linear
algebra at each level.

In physics, linear algebra is usually done over $\R$ or $\C$, but
for higher-dimensional linear algebra it is useful to start
more generally with any commutative `rig', or `ring without
negatives'.   This is a set $R$ equipped with two commutative
monoid structures, written $+$ and $\cdot$,  satisfying the
distributive law $a\cdot (b+c) = a\cdot b + a \cdot c$.  A good
example of such a thing without additive inverses is the natural
numbers (including zero), and one reason we insist on such
generality is to begin grappling with the remarkable fact that
many of the important vector spaces in physics are really
defined over the natural numbers, meaning that they contain a
canonical lattice with a basis of `positive' elements \cite{EK}.
Examples include the weight spaces of semisimple Lie algebras,
fusion algebras in conformal field theory \cite{FK}, and thanks
to the work of Kashiwara and Lusztig on canonical bases, the
finite-dimensional representations of quantum groups \cite{CP,CF}.

Linear algebra over the commutative rig $R$ can be thought of as
the study of the category $\Vect$ of `vector spaces' over $R$, by
which we simply mean $R$-modules isomorphic to $R^k$ for some
$k$.   Now $\Vect$ itself has two symmetric monoidal structures
corresponding to the direct sum $\oplus$ and tensor product
$\tensor$ of vector spaces, and the tensor product distributes
over direct sum up to a natural isomorphism satisfying certain
coherence laws.  Thus $\Vect$ is a categorical analog of a rig,
which one might call a `rig category'.  These are often called
ring categories, but there need be no additive inverses.
Precise definitions and strictification theorems for these have
been given by Laplaza \cite{Laplaza} and Kelly \cite{Kelly2}.

The analogy between the commutative rig $R$ and the symmetric rig
category $\Vect$ suggests the existence of a recursive
hierarchy of `$n$-vector spaces'.  For example, the categorical
analog of an $R$-module is a `module category' over $\Vect$.
This is a category $V$ equipped with a symmetric monoidal structure
$\oplus$ and a functor $\tensor \maps \Vect \times V \to V$
satisfying the usual conditions for a module up to natural
isomorphisms satisfying various coherence laws.
The module categories of special interest, the `2-vector spaces',
are those equivalent as module categories to the
$k$-fold Cartesian product $\Vect^k$.   A careful study of these has
been done by Kapranov and Voevodsky \cite{KV}.   The primordial
example of a 2-vector space is $\Vect$ itself, but when $R = \C$ a
more interesting example is the category of
representations of a semisimple algebra.   In general one hopes
to define a weak $(n+1)$-category  $(n+1)\Vect$ of `$(n+1)$-vector
spaces' over $R$ having as objects `module $n$-categories' over
$n\Vect$ which are $n$-equivalent, as module
$n$-categories, to $n\Vect^k$ for some $k$.   Moreover, $(n+1)\Vect$
should be a monoidal --- in fact stable --- $(n+1)$-category,
permitting one to define module $(n+1)$-categories over it and to
continue the recursive definition.   The primordial example of an
$(n+1)$-vector space should, of course, be $n\Vect$.

Much remains to be done to develop a full-fledged rigorous theory
of higher-dimensional linear algebra.   Indeed, attempting to
develop it without a well-established theory of weak
$n$-categories is rather like developing linear algebra
without solid foundations in set theory (which of course is
historically what occurred).   Nonetheless, it is already
becoming evident from examples that the sense in which
$n$-dimensional extended topological quantum field theories are
`representations' of some $n$-category is that they are
weak $n$-functors from it to $n\Vect$.

Note that while $\Vect$ is a stable 1-category, $\Hilb$ is a
stable 1-category with duals, and that the 2 levels of duality
this entails are both crucial in the definition of a unitary
TQFT.   As Freed \cite{Freed} has pointed out, one can, at least
for low $n$ so far, develop a theory of `$n$-Hilbert spaces'.
Taking our basic rig to be $\C$, recall that associated to any
vector space $V$ there is, in addition to the dual $V^\ast$, the
vector space $\overline V$ having the same additive structure but
on which $\C$ acts in a complex-conjugated manner.  A Hilbert
space is then a vector space equipped with an inner product, a
linear map
\[     \langle \cdot,\cdot \rangle \maps \overline V \tensor V \to
\C \]
that satisfies a nonnegativity condition and is also
nondegenerate, meaning that by duality it yields an isomorphism
$\overline V \cong V^\ast$.  Similarly, associated to  a 2-vector
space there is the dual $V^\ast =  \hom(V,\Vect)$, defined using
the appropriate notion of `$\hom$' for 2-vector spaces,
and also $\overline V$, the opposite category of
$V$ (i.e., having the direction of all morphisms reversed).   A
2-Hilbert space $V$ should thus be a 2-vector space equipped with
an inner product, that is, a functor
\[     \langle \cdot,\cdot
\rangle \maps \overline V \tensor V \to \Vect .\]
This should be
linear in the appropriate sense.  Nonnegativity will be automatic
because the rig category $\Vect$ has no negatives, but we should
require nondegeneracy, meaning that the inner product should
define an equivalence of $\Vect$-modules between $\overline V$
and $V^\ast$.   (That these are equivalent has already been established by
Yetter \cite{Yetter2}.)  With still further conditions giving
$V$ some of the properties of $\Hilb$, it
appears that one obtains a notion of `2-Hilbert space' such that
$2\Hilb$ is a stable 2-category with duals.  Examples should certainly
include the category of representations of a finite-dimensional
C*-algebra, such as the group algebra of a finite group
\cite{Freed}.

We hypothesize, therefore, that we can recursively define a
stable weak $n$-category with duals, `$n\Hilb$', such that the
following holds.
\vskip 1em\noindent
{\bf Extended TQFT Hypothesis, Part II.}   An
$n$-dimensional unitary extended TQFT is a weak $n$-functor,
preserving all levels of duality, from the free stable weak
$n$-category with duals on one object to $n\Hilb$.
\vskip 1em\noindent

The best evidence for this so far is the work of Freed and Quinn
on the Dijkgraaf-Witten model \cite{Freed,FQ}, of Lawrence on
extended TQFTs defined via triangulations \cite{Lawrence}, and of
Walker \cite{Walker} on Chern-Simons theory.  If this hypothesis
is true, one should be able to specify an $n$-dimensional unitary
extended TQFT simply by specifying a particular $n$-Hilbert space,
thanks to the universal property of the stable weak $n$-category with duals
on one object.  (More precisely, this should specify a weak $n$-functor
up to `equivalence' of these, a notion so far understood
only for low $n$.)  In future work we intend to describe the $n$-Hilbert
spaces giving rise to various well-known TQFTs.

\section{Quantization}

To paraphrase Nelson, quantization is a mystery, not a functor.
And indeed, while in some technical sense
we understand how quantum groups give precisely the algebraic
structures needed to construct 3-dimensional TQFTs, and how the
`quantization' of the group corresponds to the passage from
classical to quantum field theory, we do not yet know to what
extent this miracle can be generalized to higher dimensions.  The
search for algebraic structures appropriate for 4-dimensional
TQFTs is already underway, with Donaldson theory as a powerful
lure \cite{CS,CF,CKY,CY,Lawrence}.  One would like
higher-dimensional algebra to offer some guidance here, and
eventually one would very much like a comprehensive picture of
quantization for topological field theories in {\it all}
dimensions.

In its simplest guise, quantization concerns the relation between
the commutative algebras of observables in classical mechanics
and the noncommutative algebras in quantum mechanics.  On the one
hand, one can start  with a commutative algebra $A$ and try to
obtain a noncommutative algebra by `deformation quantization'.
There is in general no systematic procedure for doing this.
Nonetheless, one can study the possibilities, for example by
considering algebra structures on $A[[\hbar]]$ given by formal
power series
\[     a \star b = ab + \hbar m_1(a \tensor b) + \hbar^2 m_2(a
\tensor b) + \cdots\]
The requirement that this `star product' makes $A[[\hbar]]$
into an algebra imposes conditions on the $m_i$ which can
be studied using homological algebra \cite{BFFLS}.   Most simply,
the quantity
\[   \{a,b\} = m_1(a \tensor b) - m_1(b \tensor a),  \]
which measures the first-order deviation from commutativity of
the star product, must be a Poisson bracket on $A$, i.e.\ a
Lie bracket with
\[    \{a,bc\} = \{a,b\} c + b\{a,c\}.\]

On the other hand, there is an obvious way to get a commutative
algebra from a noncommutative algebra $A$, namely by taking its
center $Z(A)$.   Physically this corresponds to extracting the
classical part of a quantum theory, the so-called `superselection
rules'.  Note that taking the center, while a perfectly
systematic process, is {\it not} a functor, since a homomorphism
$f \maps A \to B$ need not restrict to a homomorphism from $Z(A)$
to $Z(B)$.

Figure 21 sheds some light on these ideas and their
generalizations.  In the $n = 0$ column one has, in the  $k = 1$
and $k = 2$ rows, monoids and commutative monoids respectively.
A monoid object in $\Vect$ is an algebra, while  a commutative
monoid object is a commutative algebra.  In addition to the
forgetful functor from the category of commutative algebras to
that of algebras, we have seen there are two nonfunctorial, but
still very interesting, processes relating these categories:
`taking the center' and `deformation quantization'.

These processes have analogs in other columns of Figure 21.   As
one marches down any column, one expects the last step before
stabilization to consist of imposing {\it equations} on structure
already present.  In this situation one can consider formal
deformations of the stabilized sort of structure in the category
of not-quite-stabilized structures.  For example, in the $n = 1$
column one can consider deformations of a symmetric monoidal
category in the category of braided monoidal categories.  This is
precisely where quantum groups arise!  The category of
representations of a Lie group is a symmetric monoidal category
object in $2\Vect$, while the category of representations of the
corresponding quantum group is a braided monoidal category object
in $2\Vect$.  The latter is a kind of deformation of the former
in which, for example, the braiding is given by a formal power
series
\[   B_{x,y}(a \tensor b) = b \tensor a + \hbar r_1(a
\tensor b) +  r_2(a \tensor b) + \cdots \]
The condition that $B$ be
a braiding imposes conditions on the $r_i$; for example, $r_1$
must be a solution of the `classical Yang-Baxter equations'.
For a detailed treatment of the `deformation
quantization' of symmetric monoidal categories into balanced
braided monoidal categories, see Mattes and Reshetikhin
\cite{MR}.   The tangle invariants arising this way can be
expanded as formal power series in $\hbar$, and the coefficients,
known as `Vassiliev invariants' or `invariants of finite type',
have special topological properties \cite{Barnatan,Birman}.
Their relation with the deformation quantization of commutative
algebras is clarified by the manner in which they arise in
Chern-Simons perturbation theory \cite{AMR}.

The operation of `taking the center' can also be generalized, in
a subtle and striking manner.  We can think of a $k$-tuply
monoidal $n$-category $C$ --- strict, semistrict, or weak --- as
an object in the corresponding version of $(n+k)\cat$.  Let
$Z(C)$ be the largest sub-$(n+k+1)$-category of $(n+k)\cat$
having $C$ as its only object, $1_C$ as its only morphism,
$1_{1_C}$ as its only 2-morphism, and so on, up to only one
$k$-morphism.  Thus $Z(C)$ is a $(k+1)$-tuply monoidal
$n$-category.

In what sense is $Z(C)$ the `generalized center' of $C$?
Consider first the case where $C$ is a monoid, thought of as a
category with one object.  Then $Z(C)$ is the largest
sub-2-category of $\Cat$ having $C$ as its only object, the
identity functor $1_C$ as its only morphism, and natural
transformations $\alpha\maps 1_C \doublearrow 1_C$ as 2-morphisms.  In
other words, $Z(C)$ is the commutative monoid consisting of all
natural transformations $\alpha \maps 1_C \doublearrow 1_C$.   Since
there is only one object in $C$, such a natural transformation is
simply a single morphism in $C$, and the the commutative square
condition in eq.\ (\ref{natural}) implies this morphism must
commute with all the other morphisms in $C$.  Thus $Z(C)$ is just
the center of $C$ as traditionally defined.  This also shows that
$Z$ is not a functor.

We leave it to the reader to check that if $C$ is a set, $Z(C)$
is the monoid consisting of all functions $F \maps C \to C$.
Similarly, if $C$ is a category, $Z(C)$ is the monoidal category
whose objects are functors $F \maps C \to C$ and whose morphisms
are natural transformations between such functors.  The monoidal
structure here corresponds to composition of functors.   A more
interesting example was worked out by Kapranov and Voevodsky
\cite{KV}.   Suppose that we start with a monoidal category $C$
and work in the semistrict context.    The 2-morphisms in $2\cat$
are known as `quasinatural transformations', since the square in
eq.\ (\ref{natural})  is required to commute only up to a
2-isomorphism \cite{Gray}. The 3-morphisms in $2\cat$ are known
as `modifications'.  The generalized center $Z(C)$ thus turns out
to be the braided monoidal category whose objects are
quasinatural transformations $\alpha \maps 1_C \doublearrow 1_C$
and whose morphisms are modifications between these.

It turns out that when $C$ is the monoidal category of
representations of a Hopf algebra $H$, $Z(C)$ is the braided
monoidal category of representations of a Hopf algebra called the
quantum double of $H$ \cite{CP}.  Moreover, while not themselves
quantum doubles, quantum groups are easily constructed as
quotients of quantum doubles \cite{CP}.  We thus see that that in
the $n = 1$ column of Figure 21, interesting braided monoidal
categories can be obtained either by deformation quantization of
symmetric monoidal categories, or by taking the generalized
center of monoidal categories.  We expect a more complex version
of this story to occur in the higher-$n$ columns.   For example,
in the $n = 2$ column there should be a theory of deformations of a
strongly symmetric monoidal 2-category in the category of weakly
symmetric ones.  From Figure 30 one would expect this to be
related to a Vassiliev theory for surfaces embedded in $\R^5$.
The generalized center construction should also be interesting.
For example, one can obtain braided monoidal
2-categories as the generalized centers of monoidal 2-categories
\cite{BHN}.

\subsection*{Acknowledgements}

We would like to thank Lawrence Breen, Ronald Brown, Louis Crane,
Andr\'e Joyal, Thomas Kerler, Ruth Lawrence, Stephen Sawin, Bruce
Smith, Ross Street, Dominic Verity, and David Yetter for helpful
discussions and correspondence.  We also thank Aaron Lauda for
preparing XyPic versions of all the figures in this paper.

\end{document}